\documentclass[11pt]{article}
\usepackage{amssymb,amsfonts,amsmath,latexsym,graphicx}
\newtheorem{theorem}{Theorem}[section]
\newtheorem{proposition}[theorem]{Proposition}
\newtheorem{lemma}[theorem]{Lemma}
\newtheorem{corollary}[theorem]{Corollary}
\newtheorem{definition}[theorem]{Definition}
\newtheorem{remark}[theorem]{Remark}

\newcommand{\bea}{\begin{eqnarray}}
\newcommand{\eea}{\end{eqnarray}}

\newcommand{\beq}{\begin{equation}}
\newcommand{\eeq}{\end{equation}}
\newcommand{\be}{\begin{equation}}
\newcommand{\ee}{\end{equation}}

\newcommand{\dist}{{\rm dist}}

\newcommand{\R}{{\bf R}}
\newcommand{\finedim}{{\hfill //}}
\title{
Sharp two--sided heat kernel estimates for critical
Schr\"{o}dinger \\ operators on bounded domains}

\author{\Large Stathis Filippas$^{1,4}$ \& Luisa Moschini$^{2}$
\& Achilles Tertikas$^{3,4}$  \\
                                                                           \\
        Department of Applied Mathematics$^{1}$ \\
         University of Crete,
         71409 Heraklion,  Greece \\
        filippas@tem.uoc.gr\\
                                          \\
                                          Dipartimento di Matematica ``G.
                                          Castelnuovo''$^2$,\\
 University of Rome ``La Sapienza'', 00185 Rome, Italy \\
moschini@mat.uniroma1.it \\
\\
 Department of Mathematics$^{3}$ \\
         University of Crete,
         71409 Heraklion,  Greece \\
          tertikas@math.uoc.gr\\
                                      \\
        Institute of Applied and Computational Mathematics$^4$, \\
        FORTH, 71110 Heraklion, Greece \\
    \\ }

\begin{document}
\date{}
\maketitle

\begin{abstract}
On a smooth bounded domain $\Omega\subset \R^N$ we consider the
Schr\"{o}dinger operators $-\Delta-V$, with $V$ being  either the
critical borderline potential $V(x)=(N-2)^2/4\ |x|^{-2}$ or
$V(x)=(1/4)\ \dist (x,\partial \Omega)^{-2}$, under Dirichlet
boundary conditions. In this work we  obtain sharp two-sided
estimates on the corresponding heat  kernels. To this end we
transform  the Schr\"odinger operators into suitable  degenerate
operators, for which we prove a new parabolic Harnack inequality
up to the boundary. To derive the Harnack inequality we have
established a series of new inequalities such as improved Hardy,
logarithmic Hardy Sobolev, Hardy-Moser and weighted Poincar\'e. As
a byproduct of our technique we are able to answer positively to a
conjecture of  E. B. Davies.

\smallskip

\end{abstract}

\section{Introduction and main results}

Harnack inequalities have been extremely useful in the study of
solutions of elliptic and parabolic equations, starting from the
pioneering works of De Giorgi \cite{DG}, Nash \cite{N} and Moser
\cite{Mo1}, \cite{Mo2}. They are used to prove
 H\"older continuity of solutions, strong maximum
principles, Liouville properties, as well as sharp two-sided heat
kernel estimates.  In particular, we should mention the
influential works of Aronson \cite{A} and Li and Yau \cite{LY}
where heat kernel estimates were obtained via parabolic Harnack
inequalities.

 In fact, in certain cases, the parabolic
Harnack inequality is equivalent to sharp two--sided
 heat kernel estimates.
This is the case when dealing with second order   uniformly
elliptic operators  in divergence form
 on  $\R^N$, or more generally with  weighted Laplacians on complete
Riemannian manifolds;  see the works of Fabes and Stroock
\cite{FS}, Grigoryan  \cite{G1}, and Saloff--Coste \cite{SC1}.
This equivalence has been also used in order to get sharp
two--sided
 estimates  for
Schr\"odinger operators in $\R^N$. For instance, the case of  a
potential that is regular and
 decays  like $|x|^{-2}$ at infinity
was  studied by Davies and Simon
  \cite{DS2}, where pointwise upper  bounds for the  heat
kernel were derived.  The picture was later
 completed by Grigoryan  \cite{G2} where sharp two sided estimates were
provided by means of a  parabolic Harnack inequality.
 A recent survey on heat
kernels on weighted manifolds can also  be found in \cite{G2}.

As it was shown  in
 the works of Fabes, Kenig and Serapioni  \cite{FKS}, and
Chiarenza and Serapioni   \cite{CS},  parabolic Harnack
inequalities follow after establishing Poincar\'e and Sobolev
inequalities as well as a  doubling volume growth condition.
Moreover, on complete Riemannian manifolds   parabolic Harnack
inequalities are equivalent to  Poincar\'e inequality and a
 doubling volume growth condition as explained  by Grigoryan and Saloff--Coste
in  \cite{GSC}, \cite{SC2}.

Since the  work of Baras and Goldstein \cite{BG},
 the existence or nonexistence of
solutions to  the partial differential equation \be\label{1} u_t =
\Delta u + V u, \ee with a potential $V$ involving the inverse
square of the distance function have  been widely investigated.
See  \cite{BG}, Brezis and  V\'azquez
 \cite{BV},
Cabr\'e and Martel  \cite{CM}, as well as V\'azquez and Zuazua
\cite{VZ},
 for the case  $V(x)= c|x|^{-2}$
and \cite{CM} for the  case $V(x)=c d^{-2}(x)$  on a bounded
domain $\Omega$, where $d(x) = {\rm dist}(x, \partial \Omega)$.

Concerning  the case where $V(x)= c|x|^{-2}$  with $c< (N-2)^2/4$
, sharp two--sided heat kernel estimates have been obtained in
$\R^N$, see \cite{MT1}, \cite{MT2} where the  approach of
\cite{GSC} on complete
  Riemannian manifolds  has been  used,
after a suitable transformation; see also \cite{MS} for a
different method.

On the other hand  few results are known in the case of incomplete
Riemannian manifolds, as it is for example the case of bounded
domains in $\R^N$.
 To our knowledge the only  sharp two sided estimates  in this case,
concern the standard Dirichlet Laplacian on a smooth bounded
domain $\Omega\subset \R^N$, first studied by Davies and Simon in
\cite{D1}, \cite{D2},    \cite{DS1}, and
 recently completed by Zhang \cite{Z}.  We note that in the
case of a  bounded domain,   the asymptotic of the heat kernel is
different for small time than it is for large time. In fact, for
the  heat kernel $h_{D}(t,x,y)$ of the standard Dirichlet
Laplacian and for two positive constants $C_1 \leq C_2$, we have
 for  small time
\be\label{z1} C_1 \min\left\{1, \frac{d(x)d(y)}{t} \right\}
t^{-\frac{N}{2}} e^{-C_2\frac{|x-y|^2}{t}}\le h_{D}(t,x,y) \le C_2
\min\left\{1, \frac{d(x) d(y)}{t}\right\}t^{-\frac{N}{2}}
e^{-C_1\frac{|x-y|^2}{t}}, \ee whereas for  large time
\be\label{lta} C_1 \ d(x)\  d(y) \ e^{-\lambda_1 t} \le
h_D(t,x,y)\le
 C_2\  d (x) \ d(y) \ e^{-\lambda_1 t},
\ee for all $x,y\in \Omega$; here $\lambda_1$ is the first
Dirichlet eigenvalue.

In this work, our main   interest is in obtaining
 sharp two--sided
estimates  for the  heat kernel of the
 Schr\"odinger operator  $-\Delta -V$  under  Dirichlet boundary conditions,
 on a  smooth  bounded domain
$\Omega \subset \R^N$ for the following critical borderline
potentials:
  $V(x)=
((N-2)^2/4)|x|^{-2}$  or  $V(x)=(1/4)d^{-2}(x)$.

Throughout this work $\Omega $ is a  $C^2$ bounded domain of
$\R^N$ containing the origin  and   $d(x) = {\rm dist}(x, \partial
\Omega)$. We  first consider, for  $N\geq 3$,   the  case  $V(x)=
((N-2)^2/4)|x|^{-2}$,   $x \in \Omega$
 and   we formally define the operator $K$ by
$Ku =-\Delta u-\frac{(N-2)^2}{4|x|^2}u \ ,\quad
u|_{\partial\Omega}=0$ . More precisely,  the Schr\"odinger
operator $K$ is defined in $L^2(\Omega)$ as the generator of the
symmetric form
$$\mathcal K[u_1,u_2]:=\int_{\Omega}\left(\nabla u_1 \nabla u_2 -\frac{(N-2)^2}{4 |x|^2} u_1 u_2\right) dx \ ,$$
namely, if
$$D(K):=\left\{u\in H(\Omega):-\Delta u-\frac{(N-2)^2}{4|x|^2} u \in L^2(\Omega)\right\},$$
\be\label{origin} Ku:= -\Delta u-\frac{(N-2)^2}{4|x|^2} u  \ \hbox
{ for any } u\in D(K) \ , \ee where $H(\Omega)$ denotes  the
closure of $C^\infty_0(\Omega)$ in the norm
\begin{equation}
\label{normaH} u\to ||u||_{H(\Omega)}:=\left\{\int_{\Omega}
\left(|\nabla u|^2 -\frac{(N-2)^2}{4|x|^2} u^2 \right)
 dx \right\}^{\frac{1}{2}} \ .
\end{equation}
Let us recall that $H(\Omega)\subset W^{1,q}_0(\Omega)$ for any
$1\le q<2$, due to the results in Subsection 4.1 of \cite{VZ}.

It follows, using Hardy inequality,
 that  $K$ is a nonnegative self-adjoint operator on $L^2(\Omega)$
such that for every $t>0$, $e^{-Kt}$ has an  integral kernel, that
is, $e^{-Kt}u_0(x):=\int_{\Omega} k(t,x,y)u_0(y) dy$ where
$k(t,x,y)$ is
 the heat kernel of $K$.
The first Dirichlet eigenvalue of $K$ can be defined by
\be\label{2} \lambda_1:=\inf_{0\neq \varphi \in
C^\infty_0(\Omega)} \frac{\int_{\Omega} \left(|\nabla \varphi|^2
-\frac{(N-2)^2}{4|x|^2}\varphi^2\right) dx}{\int_{\Omega}
\varphi^2 dx } \ , \ee with
 $\lambda_1
>0$,  due to  [BV]. Moreover there exists a positive function
$\varphi_1 \in H(\Omega)$ satisfying
\[
-\Delta \varphi_1-\frac{(N-2)^2}{4 |x|^2}\varphi_1=\lambda_1
\varphi_1, ~~~{\rm in}~~ \Omega,~~~ \varphi_1=0,~~  {\rm on}~~
\partial \Omega,
\]
see for example Davila and Dupaigne \cite{DD}.

We then have the following sharp two-sided heat kernel estimate on
$K$ for small time
\begin{theorem}
\label{thmorigin} Let $\Omega \subset \R^N$, $N\ge 3$, be a smooth
bounded domain containing the origin. Then there exist positive
constants $C_1, C_2$, with $C_1\le C_2$, and $T>0$ depending on
$\Omega$ such that
$$C_1 \min\Big \{(|x|+\sqrt{t})^{\frac{N-2}{2}}(|y|+\sqrt{t})^{\frac{N-2}{2}}  , \frac{d(x) d(y)}{t} \Big\}
 |xy|^{\frac{2-N}{2}} t^{-\frac{N}{2}} e^{-C_2 \frac{|x-y|^2}{t}} \le $$
$$\le k(t,x,y)\le C_2 \min\Big \{(|x|+\sqrt{t})^{\frac{N-2}{2}}(|y|+\sqrt{t})^{\frac{N-2}{2}},
\frac{d (x) d(y)}{t} \Big\}  |xy|^{\frac{2-N}{2}} t^{-\frac{N}{2}}
e^{-C_1 \frac{|x-y|^2}{t}}, $$ for all $x,y\in \Omega$ and $0<t\le
T$.
\end{theorem}
Concerning the large time asymptotic we have:

\begin{theorem}
\label{thmoriginbis} Let $\Omega \subset \R^N$, $N\ge 3$, be a
smooth bounded domain containing the origin. Then there exist two
positive constants $C_1, C_2$, with $C_1\le C_2$, such that
$$C_1 \ d(x)\  d(y) \ |xy|^{\frac{2-N}{2}} e^{-\lambda_1 t} \le k(t,x,y)\le
 C_2\  d (x) \ d(y) \ |xy|^{\frac{2-N}{2}} e^{-\lambda_1 t}, $$
for all $x,y\in \Omega$ and $t>0$ large enough; here $\lambda_1$
is defined in (\ref{2}).
\end{theorem}
To prove the above Theorem  \ref{thmoriginbis} we have shown a new
 improved Hardy inequality
which is of independent interest; see Theorem \ref{prop7}.

We next consider the case where
 the Schr\"odinger
operator $H$  has a  potential  with  critical borderline
singularity at  the boundary $Hu=-\Delta u-\frac{1}{4 d^2(x)} u\ ,
\quad u|_{\partial\Omega}=0$; here  $N\ge 2$ and  $\Omega$ is a
convex domain. More  precisely,  the Schr\"odinger operator $H$ is
defined in $L^2(\Omega)$ as the generator of the symmetric form
$$\mathcal H[u_1,u_2]:=\int_{\Omega}\left(\nabla u_1 \nabla u_2 -\frac{1}{4d^2(x)} u_1 u_2\right) dx \ ,$$
 namely if
$$D(H):=\left\{u\in W(\Omega):-\Delta u-\frac{1}{4d^2(x)} u \in L^2(\Omega)\right\},$$
\be\label{bry} Hu:= -\Delta u-\frac{1}{4d^2(x)} u  \ \hbox { for
any } u\in D(H) \ , \ee where $W(\Omega)$ denotes the closure of
$C^\infty_0(\Omega)$ in the norm
$$u\to ||u||_{W(\Omega)}:=\left\{\int_{\Omega} \left(|\nabla u|^2 -\frac{1}{4d^2(x)} u^2
 \right) dx \right\}^{\frac{1}{2}} \ .$$
Let us recall that $W(\Omega)\subset W_0^{1,q}(\Omega)$ for any
$1\le q<2$, due to Theorem B in \cite{BFT1}.

Then, due to Hardy inequality,
 $H$ is a nonnegative self-adjoint operator on $L^2(\Omega)$
such that for every $t>0$, $e^{-Ht}$ has an integral kernel, that
is, $e^{-Ht}u_0(x):=\int_{\Omega} h(t,x,y)u_0(y) dy$; here
$h(t,x,y)$ denotes
 the heat kernel of $H$.
The first Dirichlet eigenvalue of $H$ is  defined by \be\label{3}
\lambda_1:=\inf_{0\neq \varphi \in C^\infty_0(\Omega)}
\frac{\int_{\Omega} \left(|\nabla \varphi|^2
-\frac{1}{4d^2(x)}\varphi^2\right) dx}{\int_{\Omega} \varphi^2 dx
} \ . \ee It is  known  that $\lambda_1
>-\infty$ for any bounded
domain $\Omega$, and  $\lambda_1>0$ if $\Omega$ is convex, see
\cite{BM}.  Moreover there exists a positive function $\varphi_1
\in W(\Omega)$ satisfying
\[
-\Delta \varphi_1-\frac{1}{4 d^2(x)}\varphi_1=\lambda_1 \varphi_1,
~~~{\rm  in}~~ \Omega,~~~ \varphi_1=0, ~~~{\rm  on}~~~ \partial
\Omega;
\]
 see for example \cite{DD}.

We then have  the following sharp two-sided heat kernel estimate
on $H$ for small time
\begin{theorem}
\label{crit} Let $\Omega\subset \R^N$, $N \ge 2$, be a smooth
bounded and convex domain. Then there exist positive constants
$C_1, ~C_2$, with $C_1\le C_2$, and $T>0$ depending on $\Omega$
such that
$$C_1 \min\left\{1,\frac{d^{\frac{1}{2}}(x)d^{\frac{1}{2}}(y)}{t^{\frac{1}{2}}}\right\}
t^{-\frac{N}{2}} e^{-C_2\frac{|x-y|^2}{t}}\le h(t,x,y) \le C_2
\min\left\{1,\frac{d^{\frac{1}{2}}(x)
d^{\frac{1}{2}}(y)}{t^{\frac{1}{2}}}\right\} t^{-\frac{N}{2}}
e^{-C_1\frac{|x-y|^2}{t}},$$ for all $x,~y\in \Omega$ and $0<t\le
T$.
\end{theorem}
We next complement this with the large time behavior:
\begin{theorem}
\label{thmbrybis} Let $\Omega \subset \R^N$, $N\ge 2$, be a smooth
bounded and convex domain. Then there exist two positive constants
$C_1, C_2$, with $C_1\le C_2$, such that
$$C_1 \ d^{\frac{1}{2}} (x) \ d^{\frac{1}{2}}(y) \ e^{-\lambda_1 t} \le  h(t,x,y)\le
 C_2 \ d^{\frac{1}{2}} (x) \ d^{\frac{1}{2}}(y) \ e^{- \lambda_1 t}, $$
for all $x,~y \in \Omega$ and $t>0$ large enough; here $\lambda_1$
is defined in  (\ref{3}).
\end{theorem}

The two-sided estimates in Theorems   \ref{thmorigin} and
\ref{crit} are obtained as a consequence of a new parabolic
Harnack inequality
 up to the boundary, for a suitable degenerate elliptic operator.
Let  us present a  model operator  in this direction. For this we
consider classical solutions of \be\label{oper} v_t =
\frac{1}{d^\alpha(y)} div (d^\alpha(y) \nabla v), \ee (actually
solutions are considered as weak solutions, for the precise
formulation we refer to Definition \ref{sol} with $\lambda=0$
there, note that due to elliptic regularity, any solution is
smooth away from the boundary of $\Omega$).

Then, the following Harnack inequality holds true:
\begin{theorem}
\label{harnack} {\bf (Parabolic Harnack inequality
 up to the boundary)}. Let  $N\ge 2$, $\alpha\ge 1$ and $\Omega\subset \R^N$ be a smooth bounded
domain. Then, there exist positive constants $C_H$ and
$R=R(\Omega)$ such that for $x\in \Omega$, $0<r<R$ and for any
positive solution $v(y,t)$ of (\ref{oper}) in $\left\{\mathcal
B(x,r)\cap \Omega\right\}\times (0,r^2)$, the following estimate
holds true
\begin{equation} \label{solar} {\it
ess~sup}_{(y,t)\in \left\{\mathcal B(x,\frac{r}{2})\cap \Omega
\right\}\times (\frac{r^2}{4},\frac{r^2}{2})} \ \ v(y,t) \le C_H ~
{\it ess~inf}_{(y,t)\in \left\{\mathcal B(x,\frac{r}{2})\cap
\Omega \right\}\times (\frac{3}{4} r^2,r^2)} \ \ v(y,t)\
.\end{equation}
\end{theorem}
Here $\mathcal B(x,r)$ denotes roughly speaking an $N$ dimensional
cube centered at $x$ and having size $r$ (see Definition
\ref{set}). The estimate of Theorem \ref{crit} is a consequence of
Theorem \ref{harnack} and  corresponds to the extreme value
$\alpha = 1$. We refer to Theorem \ref{harnackgen} for a more
general result that leads to Theorem \ref{thmorigin}.

The existence of a uniform upper bound on the size  of the
admissible ``balls" denoted by $R(\Omega)$, is necessary, because
otherwise the nonexistence of an upper bound would imply 
two-sided heat kernel estimates that are the same for small time
and large time, which is not the case at least for $\alpha=1$, due
to Theorems \ref{crit} and \ref{thmbrybis}.

The restriction  on $\alpha$ in Theorem \ref{harnack} is sharp,
since in the  {\em weakly degenerate }  case, where $0<\alpha<1$,
even the elliptic Harnack fails. Indeed, let $\Omega:=B(0,1)$,
then $v(y):=\int^{1}_{|y|}\frac{ds}{(1-s)^\alpha s^{n-1}}$ is a
positive solution of $div(d^\alpha(y) \nabla v)=0$ for
$1/2<|y|<1$, with $v(1)=0$. The natural analogue of Theorem
\ref{harnack} in the weakly degenerate case, that is $0<\alpha<1$,
is a Harnack inequality  for the ratio of any two positive
solutions; in the elliptic case this is done by a probabilistic
approach in \cite{Ga}.

To derive heat kernel estimates we define the operator
$L:=-\frac{1}{d^\alpha(x)} div (d^\alpha(x) \nabla)$ in
$L^2(\Omega, d^\alpha(x) \ dx )$ as the generator of the symmetric
form
$$\mathcal L[v_1,v_2]:=\int_{\Omega} d^\alpha(x) \nabla v_1 \nabla v_2 \ dx \ ,$$
namely
$$D(L):=\left\{v\in H^1_0(\Omega, d^\alpha(x) \ dx ): -\frac{1}{d^\alpha(x)} div (d^\alpha(x) \nabla v)
\in L^2(\Omega, d^\alpha(x) \ dx)\right\},$$ \be\label{5} Lv:=
-\frac{1}{d^\alpha(x)} div (d^\alpha(x) \nabla v) \hbox { for any
} v\in D(L) \ , \ee
 where $H^1_0(\Omega,d^\alpha(x)
\ dx )$ denotes the closure of $C^\infty_0(\Omega)$ in the norm
\begin{equation}
\label{ehi}
 v\to ||v||_{H^1_{\alpha}}:=\left\{\int_{\Omega}
d^\alpha(x) \left(|\nabla v|^2 +v^2\right) \
dx\right\}^{\frac{1}{2}} \ .
\end{equation}
We should emphasize that for $\alpha \ge 1$, one has
$H^1_0(\Omega,d^\alpha(x) \ dx )=H^1(\Omega,d^\alpha(x) \ dx )$,
see Theorem \ref{density}.

Let us note that
 $L$ is a nonnegative self-adjoint operator on $L^2(\Omega,
d^\alpha(y) dy)$ such that for every $t>0$, $e^{-Lt}$ has a
integral kernel, that is $e^{-Lt}v_0(x):=\int_{\Omega}
l(t,x,y)v_0(y) d^\alpha(y) dy$; the existence of the heat kernel
$l(t,x,y)$   can be proved arguing as in \cite{DS1}.

Then we obtain the following sharp two-sided estimate for the heat
kernel generated by $L$.
\begin{theorem}
\label{heat} Let $\alpha\ge 1,~ N \ge 2$ and $\Omega\subset \R^N$
be a smooth bounded domain. Then there exist positive constants
$C_1, C_2$, with $C_1\le C_2$, and $T>0$ depending on $\Omega$
such that
$$C_1 \min\left\{\frac{1}{t^{\frac{\alpha}{2}}},
\frac{1}{d^{\frac{\alpha}{2}}(x)d^{\frac{\alpha}{2}}(y)}\right\}
t^{-\frac{N}{2}} e^{-C_2\frac{|x-y|^2}{t}}\le l(t,x,y) \le C_2
\min\left\{\frac{1}{t^{\frac{\alpha}{2}}},\frac{1}{d^{\frac{\alpha}{2}}(x)
d^{\frac{\alpha}{2}}(y)}\right\}t^{-\frac{N}{2}}
e^{-C_1\frac{|x-y|^2}{t}},$$ for all $x,y\in \Omega$ and $0 <t \le
T$.
\end{theorem}
%

%

So far we have considered special potentials $V$. However as we
shall see next we can obtain much more general results. For
instance we consider the operator $E:=-\Delta-V $ where the
potential $V$ is such that \be \label{h1} V(x)=V_1(x)+V_2(x)  ~ ,
\ee where \be \label{h2} |V_1(x)|\le \frac{1}{4d^2(x)} ~ , ~~~
V_2(x)\in L^p(\Omega), ~ p>\frac{N}{2} \ . \ee We also suppose
that \be \label{h3} \lambda_1:=\inf_{0\neq \varphi \in
C^\infty_0(\Omega)} \frac{\int_{\Omega} \left(|\nabla \varphi|^2
-V\varphi^2\right) dx}{\int_{\Omega} \varphi^2 dx } >0 , \ee and
that to $\lambda_1$ there corresponds a positive eigenfunction
$\varphi_1$ satisfying for all $x\in \Omega$ the following
estimate, \be \label{h4} c_1 d^{\frac{\alpha}{2}}(x) \le
\varphi_1(x)\le c_2 d^{\frac{\alpha}{2}}(x)\ , ~~ \mbox{ for some
} ~~~\alpha\ge 1\ee and for $c_1, ~ c_2$ two positive constants.

Then as before it can be shown that $E$ is a well defined
nonnegative self-adjoint operator on $L^2(\Omega)$ such that for
every $t>0$, $e^{-Et}$ has a integral kernel, that is
$e^{-Et}u_0(x):=\int_{\Omega} e(t,x,y)u_0(y) dy$. We consider
positive solutions of \be \label{hh} u_t=-Eu \ ;  \ee then our
first result reads

\begin{corollary}
\label{bbb} For $N\ge 2$, let $\Omega\subset \R^N$ be a smooth
bounded domain. Suppose that (\ref{h1}), (\ref{h2}), (\ref{h3})
and (\ref{h4}) are satisfied. Then, there exist positive constants
$C_H$ and $R=R(\Omega)$ such that for $x\in \Omega$, $0<r<R$ and
for any positive solution $u(y,t)$ of (\ref{hh}) in
$\left\{\mathcal B(x,r)\cap \Omega\right\}\times (0,r^2)$ we have
the estimate
$${\it ess~sup}_{(y,t)\in \left\{\mathcal B(x,\frac{r}{2})\cap
\Omega \right\}\times (\frac{r^2}{4},\frac{r^2}{2})} \ \ u(y,t)
d^{-\frac{\alpha}{2}}(y)
 \le C_H ~ {\it ess~inf}_{(y,t)\in \left\{\mathcal
B(x,\frac{r}{2})\cap \Omega \right\}\times (\frac{3}{4} r^2,r^2)}
\ \ u(y,t)  d^{-\frac{\alpha}{2}}(y) \ .$$
\end{corollary}

Our result in the case $\alpha=2$ is basically the local
comparison principle by Fabes, Garofalo and Salsa \cite{FGS} in
the case of Schr\"odinger operator (see Remark \ref{stab2}, that
covers the uniformly elliptic case).

As usual from the parabolic Harnack inequality the following sharp
two-sided estimate for the heat kernel $e(t,x,y)$, corresponding
to the operator $E$, can be easily deduced:
\begin{corollary}
\label{coint} For $N\ge 2$, let $\Omega\subset \R^N$ be a smooth
bounded domain.  Suppose that (\ref{h1}), (\ref{h2}), (\ref{h3})
and (\ref{h4}) are satisfied. Then there exist positive constants
$C_1, C_2$, with $C_1\le C_2$, and $T>0$ depending on $\Omega$
such that for any $x,y\in \Omega$ and $0<t\le T$
$$C_1 \min\left\{1, \frac{d^{\frac{\alpha}{2}}(x)d^{\frac{\alpha}{2}}(y)}{t^{\frac{\alpha}{2}}}\right\}
t^{-\frac{N}{2}} e^{-C_2\frac{|x-y|^2}{t}}\le e(t,x,y) \le C_2
\min\left\{1, \frac{d^{\frac{\alpha}{2}}(x)
d^{\frac{\alpha}{2}}(y)}{t^{\frac{\alpha}{2}}}\right\}t^{-\frac{N}{2}}
e^{-C_1\frac{|x-y|^2}{t}},$$ whereas for $t>T$ we have
$$C_1 \ d^{\frac{\alpha}{2}}(x)d^{\frac{\alpha}{2}}(y)
e^{-\lambda_1 t}\le e(t,x,y) \le C_2\
d^{\frac{\alpha}{2}}(x)d^{\frac{\alpha}{2}}(y) e^{-\lambda_1 t} \
.
$$
\end{corollary}
As a byproduct of our method we can answer a conjecture  by E. B.
Davies (Conjecture 7 in \cite{D2}) in the case of the
Schr\"odinger operator (see Section 4.4 for a more general case).
For this let us introduce the Green function associated to $E$,
that is \be G_E(x,y)=\int^\infty_0 e(t,x,y) dt \ , \ee then we
have
\begin{corollary}
For $N\ge 3$, let $\Omega\subset \R^N$ be a smooth bounded domain.
Suppose that (\ref{h1}), (\ref{h2}), (\ref{h3}) and (\ref{h4}) are
satisfied. Then there exist two positive constants $C_1, C_2$,
with $C_1\le C_2$, such that for any $x,y\in \Omega$
$$C_1 \min\left\{\frac{1}{|x-y|^{N-2}}, \frac{d^{\frac{\alpha}{2}}(x)d^{\frac{\alpha}{2}}(y)}{|x-y|^{N+\alpha-2}}\right\}
\le G_E(x,y) \le C_2 \min\left\{\frac{1}{|x-y|^{N-2}},
\frac{d^{\frac{\alpha}{2}}(x)
d^{\frac{\alpha}{2}}(y)}{|x-y|^{N+\alpha-2}}\right\}.$$
\end{corollary}

Davies conjecture corresponds to our result in the case
$\alpha=2$, we should note however that other values of $
\alpha\ge 1$ are possible.

The structure of the paper is as follows. In Section 2 we prove
the new parabolic Harnack inequality up to the boundary for a
doubly degenerate elliptic operator as well as the two sided heat
kernel estimates that can be deduced from it. In Section 3 we
present the proof of the above mentioned results concerning the
Schr\"odinger potential having critical singularity at the origin,
while Section 4 treats the case of the Schr\"odinger operator
having critical singularity on the boundary.

\medskip

{ \bf Acknowledgment} This work was largely done whilst the second
author was visiting the University of Crete  and FORTH  in
Heraklion, the hospitality of which is acknowledged. This research
has been partially supported by the  RTN  European network
Fronts--Singularities, HPRN-CT-2002-00274.

\setcounter{equation}{0}
\section{Parabolic Harnack inequality up to the boundary for degenerate operators}

In this section we prove a new parabolic Harnack inequality up to
the boundary for the doubly degenerate elliptic operator in
divergence form \be \label{deo}
L^\lambda_\alpha:=-\frac{1}{|x|^\lambda d^\alpha(x)} div
(|x|^\lambda d^\alpha(x) \nabla ), \ee for any $\alpha\ge 1$ and
$\lambda\in[2-N,0]$. To this end we will use the Moser iteration
technique. According to the approach presented in \cite{GSC} as
well as in \cite{CS} the key estimates one needs are the doubling
volume-growth condition (see Corollary \ref{doubling}), the local
weighted Poincar\'e inequality (see Theorem \ref{poincare}) as
well as a local weighted Moser inequality. In fact we will
establish two local weighted Moser inequalities, one will be used
near the boundary (Theorem \ref{moser}) and the other one away
from the boundary (Theorem \ref{ultima}).

Then, similarly as in \cite{GSC}, we deduce from it a sharp
two-sided heat kernel estimate. To this end a sharp volume
estimate is needed (see Lemma \ref{volume}).

\medskip

In the sequel we will use the following local representation of
the boundary of $\Omega$. There exist a finite number $m$ of
coordinate systems $(y_i',y_{iN})$, $y_i'=(y_{i1}, \cdots,
y_{iN-1})$ and the same number of functions $a_i=a_i(y_i')$
defined on the closures of the $(N-1)$ dimensional cubes
$\Delta_i:=\{y_i':|y_{ij}|\le \delta \hbox { for } j=1,\cdots,
N-1\}$, $i\in \{1,\cdots , m\}$ so that for each point $x\in
\partial \Omega$ there is at least one $i$ such that $x=(x_i',
a_i(x_i'))$. The functions $a_i$ satisfy the Lipschitz condition
on $\overline\Delta_i$ with a constant $A>0$ that is
$$|a_i(y_i')-a_i(z_i')|\le A|y_i'-z_i'|$$ for $y_i', z_i' \in
\overline \Delta_i$; moreover there exists a positive number
$\beta<1$ such that the set $B_i$ defined for any $i\in\{1,
\cdots, m\}$ by the relation $B_i=\{(y_i',y_{iN}): y_i'\in
\Delta_i, a_i(y_i')-\beta <y_{iN}<a_i(y_i')+\beta\}$ satisfy
$U_i=B_i \cap \Omega = \{(y_i',y_{iN}): y_i'\in \Delta_i,
a_i(y_i')-\beta <y_{iN}<a_i(y_i')\}$ and $\Gamma_i=B_i \cap
\partial \Omega=\{(y_i',y_{iN}): y_i'\in \Delta_i,
y_{iN}=a_i(y_i')\}$. Finally let us observe that for any $y\in
U_i$ one can make use of the following estimate on the distance
function $(1+A)^{-1} (a_i(y_i')-y_{iN})\le d(y)\le
(a_i(y_i')-y_{iN}) $  (see Corollary 4.8 in \cite{K} for details)

Let us fix from now on a constant $\gamma\in(1,2)$ and let us
define the ``balls" we will use in Moser iteration technique.
Roughly speaking they will be Euclidean balls if they stay away
from the boundary and they will be $N$ dimensional ``deformed
cubes", following the geometry of the boundary, if they are close
enough to the boundary or even if they intersect it. More
precisely we have

\begin{definition}
\label{set} (i) For any $x\in \Omega$ and for any $0<r<\beta$ we
define the ``ball" centered at $x$ and having radius $r$ as
follows. $\mathcal B(x,r)=B(x,r)$ the Euclidean ball centered at
$x$ and having radius $r$ if $d(x)\ge \gamma r$, while $\mathcal
B(x,r)=\{(y_i',y_{iN}): |y_i'-x_i'|\le r, \ a_i(y_i')-r-d(x)
<y_{iN}<a_i(y_i')+r-d(x)\}$ if $d(x)<\gamma r$, where $i\in \{1,
\cdots , m\}$ is uniquely defined by the point $\bar x \in
\partial \Omega$ such that $|\bar x-x|=d(x)$, that is by the
projection of the center $x$ onto $\partial \Omega$. (ii) We also
define by $V(x,r):=\int_{\mathcal B(x,r)\cap \Omega} |y|^\lambda
d^\alpha(y)dy$ the volume of the ``ball" centered at $x$ and
having radius $r$.
\end{definition}

We first derive a sharp volume estimate.

\begin{lemma}
\label{volume} Let $\alpha>0$, $N\ge 2$, $\lambda\in (-N,0]$ and
$\Omega\subset \R^N$ be a smooth bounded domain containing the
origin. Then there exist positive constants $c_1, c_2$ and $\beta$
such that for any $x\in \Omega$ and $0<r<\beta$, we have
$$c_1 \max\{d^\alpha (x) (|x|+r)^\lambda, r^\alpha\} r^N \le V(x,r)\le c_2
\max\{d^\alpha (x) (|x|+r)^\lambda, r^\alpha\} r^N\ .$$
\end{lemma}

To this end we make use of the following Lemma which can be proved
as in \cite{MT2}.

\begin{lemma}
\label{lele}
 Let $N\ge 2$, $\lambda\in (-N,0]$ and $\Omega\subset \R^N$
  be a smooth bounded domain containing the origin. Then there
exist two positive constants $d_1, d_2$ such that for any $x\in
\Omega$, we have
\begin{equation}\label{cita1} d_1 \ r^N(|x|+r)^\lambda \le
\int_{B(x,r)}|y|^\lambda dy \le d_2 \ r^N(|x|+r)^\lambda\
.\end{equation}
\end{lemma}

Let us accept (\ref{cita1}) at the moment and let us prove the
sharp volume estimate.

{\em Proof of Lemma \ref{volume}: } Let us first consider the case
where $d(x)\ge \gamma r$. Then $\mathcal B(x,r)=B(x,r)\subset
\Omega$. Due to the fact that any $y\in B(x,r)$ satisfies
\beq\label{bee} \left(\frac{\gamma-1}{\gamma}\right) d(x)\le
d(x)-r\le d(y)\le d(x)+r \le \left(\frac{\gamma+1}{\gamma}\right)
d(x)\eeq the claim easily follows making use of Lemma \ref{lele}
with $c_2\ge d_2 \left(\frac{\gamma+1}{\gamma}\right)^\alpha$ and
$c_1\le d_1 \left(\frac{\gamma-1}{\gamma}\right)^\alpha$.

Let us now consider the case where $d(x)<\gamma r$ and let us
denote by $L_1,L_2$ two positive constants such that
$B(0,L_1)\subset \Omega\subset B(0,L_2)$ (note that they exist by
assumption on $\Omega$). Then any $y\in \mathcal B(x,r)\cap
\Omega$ satisfies the following estimate $L_1 -(\gamma+1) \beta\le
|y|\le L_2$. Indeed, if on the contrary $|y|<L_1-(\gamma+1)
\beta$, then by definition of $L_1$ we would have $d(y)>(\gamma+1)
\beta$, and this contradicts our assumption. In fact one obviously
has $d(y)\le d(x)+r< (\gamma+1)r<(\gamma+1)\beta$ and it is not
restrictive to suppose from the beginning that the parameter
$\beta$ in the local representation of the boundary of $\Omega$
satisfies $\beta<L_1(\gamma+1)^{-1}$. As a consequence we have:
\beq \label{ohi} \forall \ y\in \mathcal B(x,r)\cap \Omega~ ~ ~
d_3 \le |y|^\lambda\le d_4 \ .\eeq here $d_3:=L_2^\lambda$ and
$d_4:=(L_1 -(\gamma+1) \beta)^\lambda$.

Then for some $i\in \{1,\cdots,m\}$, we have
$$V(x,r)=\int_{\mathcal B(x,r)\cap \Omega} |y|^\lambda d^\alpha (y) dy\sim
\int_{|y_i'-x_i'|\le
r}\int^{\min\{a_i(y_i'),a_i(y_i')+r-d(x)\}}_{a_i(y_i')-r-d(x)}
|y|^\lambda (a_i(y_i')-y_{iN})^\alpha dy_{iN} dy_i' \ .$$ From now
on we omit the subscript $i$ for convenience.
 Indeed we have
$$V(x,r)\le \int_{|y'-x'|\le r}\int^{\min\{a(y'),a(y')+r-d(x)\}}_{a(y')-r-d(x)}
|y|^\lambda (a(y')-y_{N})^\alpha dy_N dy'\le d_4 (\gamma+1) r
(d(x)+r)^{\alpha} \int_{|y'-x'|\le r} dy'\le $$
$$\le d_4 (\gamma+1)^{\alpha+1} r^{\alpha+1+N-1}
  \omega_{N-1}=d_4 (\gamma+1)^{\alpha+1} r^{\alpha+N} \omega_{N-1}\ .$$
On the other hand
$$V(x,r)\ge  (1+A)^{-\alpha} \int_{|y'-x'|\le r}\int^{a(y')-r(\gamma-1)}_{a(y')-r-d(x)} |y|^\lambda (a(y')-y_{N})^\alpha
  dy_N
dy' \ge $$ $$\ge d_3 (1+A)^{-\alpha}\int_{|y'-x'|\le
r}\int^{a(y')-r(\gamma-1)}_{a(y')-r-d(x)} (r(\gamma-1))^\alpha
dy_N dy'\ge
$$
$$\ge d_3 (1+A)^{-\alpha}r^{\alpha+1} (2-\gamma) (\gamma-1)^{\alpha} \int_{|y'-x'|\le r}dy'=
d_3 (1+A)^{-\alpha}r^{\alpha+1+N-1}
(2-\gamma)(\gamma-1)^{\alpha}\omega_{N-1}=$$ $$=d_3
(1+A)^{-\alpha}r^{\alpha+N} (2-\gamma)
(\gamma-1)^{\alpha}\omega_{N-1}\ . $$ Here $\omega_{N}$ denotes
the standard volume of the Euclidean unit ball in $\R^N$. Thus the
result follows with $c_1:=\min\left\{d_1
\left(\frac{\gamma-1}{\gamma}\right)^{\alpha},
d_3(1+A)^{-\alpha}(2-\gamma)(\gamma-1)^{\alpha}\omega_{N-1}\right\}$
and $c_2:=\max
\left\{d_2\left(\frac{\gamma+1}{\gamma}\right)^{\alpha},
d_4(\gamma+1)^{\alpha+1}\omega_{N-1}\right\}$.

\finedim

Let us now prove estimate (\ref{cita1}), which is taken from
\cite{MT2}, we give here the details for the convenience of the
reader

{\em Proof of Lemma \ref{lele}: } $(i)$ Observe that
$$
 \int_{B(x,r)}|y|^\lambda dy =r^{\lambda+N} \int_{B(0,1)}|w+z|^\lambda dz \ ,
 $$
where $w:=\frac{x}{r}$. Then (\ref{cita1}) reads
 \be\label{stimanum}
 d_1 (|w|+1)^\lambda \le \int_{B(0,1)}|w+z|^\lambda dz \le d_2
(|w|+1)^\lambda \ . \ee
  Since $|z|\le 1$,
 there holds
$$
 |w|-1\le |w|-|z|  \le |w+z| \le |w| +|z| \le |w|+1\ ,
 $$
  whence
 \be\label{stin}
 \omega_N \Big(|w|+1 \Big)^\lambda \le \int_{B(0,1)}|w+z|^\lambda dz \le
 \omega_N \Big(|w|-1 \Big)^\lambda \ ;
 \ee
Comparing (\ref {stimanum}) and (\ref{stin}), we immediately
obtain the lower bound in (\ref{cita1})  with $d_1:= \omega_N$.
(ii) To prove the upper bound in (\ref{cita1}), observe that three
cases are possible: (a) $r\le \frac{|x|}{2}$, (b)
$\frac{|x|}{2}<r<3|x|$ and (c) $r\ge 3|x|$. In case (a) the claim
follows from the right-hand inequality in (\ref{stin}), if we
exhibit $d_2>0$ such that
 $$
 \omega_N(|w|-1)^\lambda \le d_2 (|w|+1)^\lambda
 $$
 for any $|w|\ge 2$.
It is easily seen that the function
$F(t):=\left(\frac{t+1}{t-1}\right)^{|\lambda|} \; (t>1)$ is
decreasing, thus $F(t)\le F(2)=\left(\frac{1}{3}\right)^{\lambda}$
for any $t\ge 2$. This proves the claim in this case for any
$d_2\ge \omega_N \frac{1}{3^{\lambda}}$.

To deal with cases $(b)-(c)$, observe first that
$$
\int_{B(0,r)}|y|^\lambda dy =\frac{\omega_N}{\lambda+N}r^{\lambda
+N} \ .
$$
In case $(b)$ we have \beq\label{incl1} B(x,r)\subseteq
B\left(0,4|x|\right) \; , \eeq whence \beq\label{doux1}
\int_{B(x,r)}|y|^\lambda dy \le \frac{\omega_N}{\lambda+N}(4
|x|)^{\lambda +N} \le \frac{\omega_N}{\lambda+N}(8 r)^{\lambda +N}
\ . \eeq In  case $(c)$ there holds \beq\label{incl2}
B(x,r)\subseteq B\left(0,\frac{4 r}{3}\right) \  , \eeq thus
\beq\label{doux2} \int_{B(x,r)}|y|^\lambda dy \le
\frac{\omega_N}{\lambda+N}\left(\frac{4r}{3}\right)^{\lambda +N} \
. \eeq Since $$r^{\lambda+N}\le \frac{r^N (|x|+r)^\lambda
}{3^\lambda} \ ~ ~ \hbox { for } |x|<2r$$ we have
$$\int_{B(x,r)}|y|^\lambda dy \le
\frac{\omega_N}{\lambda+N}\frac{8^{\lambda+N}}{3^\lambda} r^N
(|x|+r)^\lambda $$ in cases (b)-(c). Hence the conclusion follows
with $d_2:= \frac{\omega_N}{\lambda
+N}\frac{8^{\lambda+N}}{3^\lambda}$.

\finedim

>From Lemma \ref{volume} one can easily deduce the doubling
property which reads as follows:

\begin{corollary}
\label{doubling} Let $\alpha>0$, $N\ge 2$, $\lambda\in (-N,0]$ and
$\Omega\subset \R^N$ be a smooth bounded domain containing the
origin. Then
 there exist positive constants $C_D$ and $\beta$ such that for any $x\in \Omega$ and
$0<r<\beta$, we have
$$V(x,2r)\le C_D V(x,r).$$
\end{corollary}

Let us state now the local Poincar\'e inequality.

\begin{theorem}{\bf (Local weighted Poincar\'e inequality)}
\label{poincare} Let $\alpha>0$, $N\ge 2$, $\lambda\in (-N,0]$ and
$\Omega\subset \R^N$ be a smooth bounded domain containing the
origin. Then there exist positive constants $C_P$ and $\beta$ such
that for any $x\in \Omega$ and $0<r<\beta$, we have
\begin{equation} \label{poipoi}\inf_{\xi \in \R}\int_{\mathcal
B(x,r)\cap \Omega}|y|^\lambda d^\alpha(y)|f(y)-\xi|^2  dy\le C_P \
r^2 \int_{\mathcal B(x,r)\cap \Omega} |y|^\lambda d^\alpha(y)
|\nabla f|^2 dy \ , ~ ~ ~ \forall \ f \in C^1(\overline {\mathcal
B(x,r)\cap \Omega}) \ .
\end{equation}
\end{theorem}

{\em Proof: } Let us first consider the case where $d(x)\ge \gamma
r$. Then $\mathcal B(x,r)=B(x,r)\subset\Omega$. Due to (\ref{bee})
the claim corresponds to Theorem 3.1 in \cite{MT2}. We give here
the details for the convenience of the reader.

(i) As a consequence of the compact embedding of the space
$H^1(B(0,1), |y|^\lambda dy)$ into $L^2(B(0,1), |y|^\lambda dy)$
(e.g. see \cite{KO}) we have that
$$\int_{B(0,1)} |f-\hat f|^2 |y|^\lambda dy \le C
\int_{B(0,1)}|\nabla f|^2 |y|^\lambda dy \ , ~ ~ ~ \forall \ f\in
C^1(\overline{B(0,1)}) \ ;$$ here $\hat f:=\left(\int_{B(0,1)}
f(y) |y|^\lambda dy\right)\left(\int_{B(0,1)} |y|^\lambda
dy\right)^{-1}$. Then by scaling (\ref{poipoi}) follows when
$x=0$.

(ii) Let us now consider the case $|x|\ge 2r$ and let us define
$\bar f:=\omega_N^{-1}\int_{B(0,1)} f(x+rz) dz$. Then as in the
proof of Lemma \ref{lele} we have $(|x|+r)^\lambda \le
|x+rz|^\lambda \le \frac{1}{3^\lambda}(|x|+r)^\lambda$. Hence
$$\int_{B(x,r)}
|f-\bar f|^2 |y|^\lambda dy = r^N \int_{B(0,1)} |f(x+rz)-\bar f|^2
|x+rz|^\lambda dz \le
\frac{r^N}{3^\lambda}(|x|+r)^\lambda\int_{B(0,1)} |f(x+rz)-\bar
f|^2 dz \le $$ $$\le C r^2
\frac{r^N}{3^\lambda}(|x|+r)^\lambda\int_{B(0,1)} |\nabla f|^2
(x+rz) dz \le C r^2 \frac{r^N}{3^\lambda} \int_{B(0,1)} |\nabla
f|^2 |x+rz|^\lambda dz = C r^2 \frac{1}{3^\lambda} \int_{B(x,r)}
|\nabla f|^2 |y|^\lambda dy \ .$$ Then (\ref{poipoi}) follows when
$|x|\ge 2r$.

(iii) For a general $x\in \Omega$ two cases are possible (a) $0\le
|x|<\frac{r}{4}$; (b) $|x|\ge \frac{r}{4}$. In case (a) there
holds
$$B\left (x,\frac{r}{8}\right)\subseteq
B\left(0,\frac{r}{2}\right)\subseteq B(x,r) \ ,$$ thus from (i) we
have
$$\inf_{\xi\in \R} \int_{B\left(x,\frac{r}{8}\right)}
|f(y)-\xi|^2|y|^\lambda dy \le\inf_{\xi\in \R}
\int_{B\left(0,\frac{r}{2}\right)} |f(y)-\xi|^2|y|^\lambda dy \le
\frac{C}{4}r^2 \int_{B\left(0,\frac{r}{2}\right)} |\nabla
f|^2|y|^\lambda dy \le C_P r^2 \int_{B(x,r)} |\nabla f|^2
|y|^\lambda dy \ ,
$$
This proves (\ref{poipoi}) in case (a) since using a Whitney type
covering and arguing as in \cite{SC2} the integration set of the
left hand side which from above is $B\left(x,\frac{r}{8}\right)$
can be increased as to cover all $B(x,r)$.

In case (b) there holds $|x|\ge 2\left(\frac{r}{8}\right)$; hence
from (ii)
$$\inf_{\xi\in \R} \int_{B\left(x,\frac{r}{8}\right)}
|f(y)-\xi|^2|y|^\lambda dy \le \frac{C}{64}r^2
\int_{B\left(x,\frac{r}{8}\right)} |\nabla f|^2|y|^\lambda dy \le
C_P r^2 \int_{B(x,r)} |\nabla f|^2 |y|^\lambda dy \ .
$$
This completes the proof in the case $d(x)\ge \gamma r$.

Let us now consider the case where $d(x)<\gamma r$. Then for some
$i\in \{1,\cdots,m\}$ we have $$\int_{\mathcal B(x,r)\cap
\Omega}d^\alpha (y)|f(y)-\xi|^2  dy\le  \int_{|y_i'-x_i'|\le
r}\int^{\min\{a_i(y_i'), a_i(y_i')+r-d(x)\}}_{a_i(y_i')-r-d(x)}
|f(y_i',y_{iN})-\xi|^2 |y|^\lambda (a_i(y_i')-y_{iN})^\alpha
dy_{iN} dy_i'$$ From now on we omit the subscript $i$ for
convenience. Let us perform then the following change of variables
$(y',y_N)\to (y',z_N:=a(y')-y_N)$ and make use of (\ref{ohi});
thus the above integral is less or equal then:
$$d_4 \int_{|y'-x'|\le r}\int_{\max\{0,d(x)-r\}}^{r+d(x)} |f(y',a(y')-z_{N})-\xi|^2 z_N^\alpha
\ dz_{N} dy'=$$ $$= d_4\int_{\max\{0,d(x)-r\}}^{r+d(x)} z_N^\alpha
\left(\int_{|y'-x'|\le r}|f(y',a(y')-z_{N})-\xi|^2
dy'\right)dz_{N}\le $$
$$\le C d_4 r^2\int_{\max\{0,d(x)-r\}}^{r+d(x)} z_N^\alpha
\left(\int_{|y'-x'|\le r}\Big|\frac{\partial f}{\partial
y'}+\frac{\partial f}{\partial y_n} \frac{a(y')}{\partial
y'}\Big|^2  dy'\right)dz_{N} \le$$ $$\le C d_4r^2
\int_{\max\{0,d(x)-r\}}^{r+d(x)} \int_{|y'-x'|\le r}|\nabla
f|^2(y',a(y')- z_N) z_N^\alpha\  dy'dz_{N}\le C
\frac{d_4}{d_3}(1+A)^\alpha r^2 \int_{\mathcal B(x,r)\cap \Omega}
|y|^\lambda d^\alpha (y)|\nabla f|^2  dy \ .$$

In the above argument we made the following choices
$$\xi=\xi(z_N):=\left(\int_{|y'-x'|\le r} f(y', a(y')-z_N) dy'\right)r^{-N+1} \omega_{N-1}^{-1}
= \omega_{N-1}^{-1}\int_{|z'-x'|\le 1} f(rz',a(rz')-z_N) dz'\ ,$$
and $C$ being the Euclidean Poincar\'e constant on the $N-1$
dimensional Euclidean ball of radius one.

Since for any $\bar{\xi} \in \R$, $|f-\bar{\xi}|^2\le 2
|f-\xi(z_N)|^2 +2|\xi(z_N)-\bar{\xi}|^2$ in order to prove
(\ref{poipoi}) in this case, it only remains to estimate the
following term
$$\int_{|y'-x'|\le r}\int_{\max\{0,d(x)-r\}}^{r+d(x)} |\xi(z_{N})-\bar{\xi}|^2 z_N^\alpha \ dz_{N}
dy'= $$
$$=\left(\int_{|y'-x'|\le r}dy'\right)\Big[|\xi(z_{N})-\bar{\xi}|^2
 \frac{z_N^{\alpha+1}-\max\{0,d(x)-r\}^{\alpha+1}}{\alpha+1}\Big|_{\max\{0,d(x)-r\}}^{r+d(x)}
 -$$
 $$-
\frac{2}{\alpha+1} \int_{\max\{0,d(x)-r\}}^{r+d(x)}
(\xi(z_{N})-\bar \xi)\frac{\partial \xi(z_N)}{\partial z_N} \left(
z_N^{\alpha+1}-  \max\{0,d(x)-r\}^{\alpha+1}\right) dz_{N} \Big]\
.$$ Thus, choosing $\bar \xi:=\xi(r+d(x))$ above, we obtain by
H\"older inequality
$$\int_{|y'-x'|\le r}\int_{\max\{0,d(x)-r\}}^{r+d(x)} |\xi(z_{N})-\bar \xi|^2 z_N^\alpha \ dz_{N} dy'\le $$
$$\le \frac{2}{\alpha+1}
\left(\int_{|y'-x'|\le r}\int_{\max\{0,d(x)-r\}}^{r+d(x)}
|\xi(z_{N})-\bar \xi|^2 z_N^\alpha \ dz_{N}
dy'\right)^{\frac{1}{2}}\times$$ $$\times \left(\int_{|y'-x'|\le
r}\int_{\max\{0,d(x)-r\}}^{r+d(x)} \Big| \frac{\partial
\xi(z_N)}{\partial z_N} \Big|^2 z_N^{\alpha}
\left(\frac{z_N^{\alpha+1}-\max\{0,d(x)-r\}^{\alpha+1}}{z_N^{\alpha}}\right)^2
dz_{N} dy'\right)^{\frac{1}{2}}\ .
$$
Since $$\Big| \frac{\partial \xi(z_N)}{\partial z_N} \Big|\le
\omega_{N-1}^{-1} \int_{|z'-x'|\le 1} \Big| \frac{\partial
f}{\partial y_N}\Big|(rz',a(rz')-z_N)dz' =
\omega_{N-1}^{-1}r^{-N+1}\int_{|y'-x'|\le r} \Big| \frac{\partial
f}{\partial y_N}\Big|(y',a(y')-z_N)dy'\ ,$$ hence
$$\Big|
\frac{\partial \xi(z_N)}{\partial z_N} \Big|^2\le
\omega_{N-1}^{-1}r^{1-N}\int_{|y'-x'|\le r} \Big| \frac{\partial
f}{\partial y_N}\Big|^2(y',a(y')-z_N)dy'\ .$$
 Thus, since $d(x)<\gamma r$, we obtain (\ref{poipoi}) with constant
$C_P:=2\frac{d_4}{d_3}(1+A)^\alpha
\left(C+\frac{4(\gamma+1)^2}{(\alpha+1)^2}\right)$.

\finedim

\medskip
We next prove the following
 local weighted Moser inequality:

\smallskip
\begin{theorem}
{\bf (Local weighted Moser inequality)} \label{moser} Let
$\alpha>0$, $N\ge 2$ and $\Omega\subset \R^N$ be a smooth bounded
domain. Then there exist positive constants $C_M$ and
$R=R(\alpha,\Omega)$ such that for any $\nu\ge N+\alpha$, $x\in
\Omega$, $0<r<R$ and $f \in C^\infty_0(\mathcal B(x,r))$ we have
$$\int_{\mathcal B(x,r)\cap \Omega}d^\alpha(y)|f(y)|^{2\left(1+\frac{2}{\nu}\right)}  dy\le C_M r^2
V(x,r)^{-\frac{2}{\nu}} \left(\int_{\mathcal B(x,r)\cap
\Omega}d^\alpha(y) |\nabla f|^2 dy \right) \left(\int_{\mathcal
B(x,r)\cap \Omega} d^\alpha(y) |f|^2 dy \right)^\frac{2}{\nu} .
$$
\end{theorem}
{\em Proof: } Let us first consider the
case where $d(x)\ge \gamma r$. By the standard Moser inequality,
there exists a positive constant $C$ such that for any $x\in
\Omega$ and any $\nu\ge N$ if $N\ge 3$ or any $\nu>2$ if $N=2$,
the following holds true
$$\int_{B(x,r)}|f(y)|^{2\left(1+\frac{2}{\nu}\right)} dy\le C r^2
r^{-\frac{2N}{\nu}} \left(\int_{B(x,r)} |\nabla f|^2 dy \right)
\left(\int_{B(x,r)} |f|^2 dy \right)^\frac{2}{\nu} \ \ , ~ ~ ~
\forall \ f \in C^\infty_0(B(x,r))$$ (see for example Section
2.1.3 in \cite{SC2}). Thus we have
$$\int_{B(x,r)}d^\alpha(y) |f(y)|^{2\left(1+\frac{2}{\nu}\right)}  dy\le (d(x)+r)^\alpha C r^2
r^{-\frac{2N}{\nu}} \left(\int_{B(x,r)} |\nabla f|^2 dy \right)
\left(\int_{B(x,r)} |f|^2 dy \right)^\frac{2}{\nu}\le $$
$$\le Cr^2 \left(\frac{d(x)+r}{d(x)-r}\right)^\alpha (r^N (d(x)-r)^\alpha)^{-\frac{2}{\nu}}
\left(\int_{B(x,r)} d^\alpha(y) |\nabla f|^2 dy \right)
\left(\int_{B(x,r)} d^\alpha(y) |f|^2  dy \right)^\frac{2}{\nu}
\le
$$
$$\le C_M r^2 V(x,r)^{-\frac{2}{\nu}} \left(\int_{B(x,r)} d^\alpha(y) |\nabla f|^2 dy \right)
\left(\int_{B(x,r)} d^\alpha(y) |f|^2 dy \right)^\frac{2}{\nu}\
,$$ where $C_M:=C \left(1+\frac{2}{\gamma-1}\right)^\alpha
\left(\frac{\gamma}{\gamma-1}\right)^\frac{2\alpha}{\nu}
c_2^\frac{2}{\nu}$ and $c_2$ is the constant appearing in the
volume estimate in Lemma \ref{volume} when $\lambda=0$.

Let us now consider the case where $d(x)<\gamma r$. Then we claim
the following local weighted Sobolev inequality: there exist
positive constants $C_S$ and $R=R(\alpha,\Omega)$ such that for
any $x\in \Omega$, $0<r<R$, satisfying $d(x)<\gamma r$, and any
$f\in C^\infty_0(\mathcal B(x,r))$, we have
\begin{equation} \label{ls} \left(\int_{\mathcal B(x,r)\cap
\Omega}d^\alpha(y) |f(y)|^{\frac{2(N+\alpha)}{N+\alpha-2}}
dy\right)^{\frac{N+\alpha-2}{N+\alpha}}\le C_S \int_{\mathcal
B(x,r)\cap \Omega} d^\alpha(y)|\nabla f|^2  dy \ .\end{equation}

If we accept (\ref{ls}), then the result follows, with $C_M=C_S$
by means of H\"older inequality, in fact we have
$$\int_{\mathcal B(x,r)\cap \Omega}
d^\alpha (y) f^{2\left(1+\frac{2}{N+\alpha}\right)} dy \le
\left(\int_{\mathcal B(x,r)\cap \Omega} d^\alpha
(y)f^{\frac{2(N+\alpha)}{N+\alpha-2}}
dy\right)^{\frac{N+\alpha-2}{N+\alpha}}\left(\int_{\mathcal
B(x,r)\cap \Omega} d^\alpha (y)f^{2}
dy\right)^{\frac{2}{N+\alpha}}$$ as well as for any $\nu>
N+\alpha$
$$\left(\int_{\mathcal B(x,r)\cap \Omega}
d^\alpha (y)f^{2\left(1+\frac{2}{\nu}\right)} dy\right)
\left(\int_{\mathcal B(x,r)\cap \Omega} d^\alpha (y)f^{2}
dy\right)^{\frac{2(\nu-N-\alpha)}{\nu(N+\alpha)}} \le
\left(\int_{\mathcal B(x,r)\cap \Omega} d^\alpha
(y)f^{2\left(1+\frac{2}{N+\alpha}\right)} dy\right)
V(x,r)^{\frac{2(\nu-N-\alpha)}{\nu(N+\alpha)}} \ .$$ In the sequel
we will give the proof of (\ref{ls}). We will follow closely the
argument of \cite{FMT2}. If $V\subset \R^N$ is any bounded domain
and $u\in C^\infty(\overline V)$ then it is well known that
$$S_N ||u||_{L^{\frac{N}{N-1}}(V)} \le ||\nabla u||_{L^1(V)} +||u||_{L^1(\partial V)} \ ,$$
where
$S_N:=N\pi^{\frac{1}{2}}\left[\Gamma(1+\frac{N}{2})\right]^{-\frac{1}{2}}$
(see p. 189 in \cite{M}).  Let us fix from now on that
$V:=\mathcal B(x,r)\cap \Omega$, and let us apply the above
inequality to $u:=d^a f$,  for any $f\in C^\infty_0(\mathcal
B(x,r))$ and any $a>0$. Thus we get
$$S_N ||d^a f||_{L^{\frac{N}{N-1}}(V)} \le \int_{V}\left(|\nabla f|d^a +a d^{a-1}|\nabla d||f|\right)
dy \ .$$ Let us remark at this point that boundary terms on
$\partial \Omega$ are zero due to the presence of the weight
$d^a$, $a>0$. To estimate the last term of the right hand side, we
will make use of an integration by parts, noting that $\nabla
d\cdot \nabla d=1$ a.e.; that is we have:
$$a\int_V d^{a-1} |f| dy =a\int_V \nabla d \cdot \nabla d \ d^{a-1} |f| dy=
\int_{V}\nabla d^a \cdot \nabla d |f| dy=
$$
$$=-\int_{V} d^a \Delta d |f| dy -\int_{V} d^a \nabla d\cdot \nabla |f| dy +
\int_{\partial V} d^a  \nabla d \cdot \nu \ |f|\ dS .$$

Under our smoothness assumption on $\Omega$ we have that $|d\Delta
d|\le c_0 \delta$ in $\Omega_\delta$ for $\delta$ small, say
$0<\delta\le \delta_0$,  and for some positive constant $c_0$
independent of $\delta$ ($\delta_0$, $c_0$ depending on $\Omega$).
Now, if $d(x)+r<\delta$, that is if $r<\frac{\delta}{\gamma+1}$,
we have that $V \subset \Omega_\delta $ and it follows that
$$a\int_V d^{a-1} |f| dy\le c_0 \delta \int_V d^{a-1} |f| dy + \int_V d^{a} |\nabla f| dy \ ,$$
hence
\begin{equation}
\label{ena} \int_V d^{a-1} |f| dy\le (a-c_0\delta)^{-1}\int_V
d^{a} |\nabla f| dy \ .
\end{equation}
Consequently for any $r\in (0,R(a,\Omega))$,
$R(a,\Omega):=\frac{1}{\gamma+1}\min\{\delta_0,\frac{a}{c_0}\}$
and any $a>0$ the following inequality is true
\begin{equation}
\label{deca} S_N ||d^a f||_{L^{\frac{N}{N-1}}(V)} \le
\left(\frac{a}{a-c_0\delta_0} +1\right) \int_{V}d^a|\nabla f| dy \
.\end{equation}

To proceed we will use the following interpolation inequality (cf.
Lemma 4.1 of \cite{FMT2}).

$$||d^bf||_{L^q(V)}\le \frac{N(q-1)}{q} ||d^a f||_{L^{\frac{N}{N-1}}(V)} +
\frac{q-N(q-1)}{q}||d^{a-1} f||_{L^1(V)} \ ,$$

\begin{equation} \label{in0}
\forall \ 1<q\le \frac{N}{N-1} \ ,~ ~ ~  b:=a-1+\frac{q-1}{q}N \
,~ ~ ~  a>0\ . \end{equation}

>From (\ref{ena}) and (\ref{deca}), we get for any $a,b,q$ as above
the following inequality
\begin{equation}
\label{in1} ||d^bf||_{L^q(V)}\le C_1||d^a \nabla f||_{L^1(V)} \ ,
\end{equation}
where $C_1:=\frac{N(q-1)}{q}\frac{1}{S_N}
\left(\frac{a}{a-c_0\delta_0}+1\right)+\frac{q-N(q-1)}{q}\left(\frac{1}{a-c_0\delta_0}\right)$.

Let us now apply inequality (\ref{in1}) to $|f|^s$ instead of $f$,
for $s:=\frac{Q}{2}+1$, $q:=\frac{Q}{s}$, $b:=Bs$.  Due to
(\ref{in0}) we have $ a=b+1-\frac{q-1}{q}N=\frac{BQ}{2}+A$, where
$A:=B+1 -\frac{Q-2}{2Q}N$. In this way we obtain
$$\left(\int_{V} d^{BQ} |f|^Q dy \right)^{\left(\frac{Q}{2}+1\right)\frac{1}{Q}} \le
C_1  \left(\frac{Q}{2}+1\right) \left(\int_{V} d^{\frac{BQ}{2}+A}
|f|^{\frac{Q}{2}} |\nabla f| dy \right) \le$$ $$\le
C_2\left(\int_{V} d^{BQ} |f|^Q dy\right)^{\frac{1}{2}}
\left(\int_{V} d^{2A} |\nabla f|^2 dy \right)^{\frac{1}{2}} \ ;$$
where $C_2:=C_1\left(\frac{Q}{2}+1\right)$.

After simplifying we see that we have proved the following: there
exists $R=R(\frac{BQ}{2}+A,\Omega)$ such that for all $0<r<R$ and
all $x\in \Omega$ with $d(x)<\gamma r$, there holds
$$\left(\int_{\mathcal B(x,r)\cap \Omega} d^{BQ} |f|^{Q} dy \right)^{\frac{2}{Q}}
\le C \int_{\mathcal B(x,r)\cap \Omega}d^{2A}(y) |\nabla f|^2 dy \
,$$ for any $N\ge 2$ and any $f \in C^\infty_0(\mathcal B(x,r))$
under the following conditions $A:=B+1-\frac{Q-2}{2Q}N$,
$\frac{BQ}{2}+A>0$, $2<Q<\infty$ if $N=2$, $2<Q\le \frac{2N}{N-2}$
if $N\ge 3$; here $C_3=C_2^2=C_3(N,Q,B,c_0,\delta_0)$.

Taking $A=\frac{\alpha}{2}$, $Q:=\frac{2(N+\alpha)}{N+\alpha-2}$
and $B:= \frac{\alpha}{Q}$
we deduce the local weighted Sobolev inequality (\ref{ls}) with
$C_S=C_S(N,\alpha,c_0,\delta_0)$ and this completes the proof of
Theorem \ref{moser}.

\finedim

\begin{remark}
Note that the upper bound for the length of the ``balls" in the
local weighted Moser inequality, denoted by $R(\alpha,\Omega)$,
goes to zero as $\alpha$ tends to zero.
\end{remark}

\begin{remark}
Let us note that when $N=1$, the corresponding analogue of the
local weighted Sobolev inequality  (\ref{ls}) when $\Omega=(-1,1)$
is the following one
$$\left(\int^{\min\{1,x+r\}}_{\max\{-1,x-r\}} (1-|y|)^\alpha |f|^q(y) dy
\right)^{\frac{1}{q}}\le C_S\
r^{\frac{\alpha+1}{q}+\frac{1-\alpha}{2}}\left(\int^{\min\{1,x+r\}}_{\max\{-1,x-r\}}
(1-|y|)^\alpha |f'|^2(y) dy \right)^{\frac{1}{2}}\ ,$$ for any
$f\in C^\infty_0(x-r,x+r)$, and any $q>2$ if $0<\alpha\le 1$ and
$2<q\le \frac{2(\alpha+1)}{\alpha-1}$ if $\alpha>1$. Consequently
Theorem \ref{harnack} as well as its consequences can be also
stated for $N=1$.
\end{remark}

>From the results within this subsection, we will now deduce a new
parabolic Harnack inequality up to the boundary for the doubly
degenerate elliptic operator $L^\lambda_\alpha$  defined in
(\ref{deo}). To this end let us first make precise the notion of a
weak solution
\begin{definition}
\label{sol} By a solution $v(y,t)$ to $v_t =-L^\lambda_\alpha v$
in $Q:=\{\mathcal B(x,r)\cap \Omega \}\times (0,r^2)$, we mean a
function $v\in C^1((0,r^2); L^2(\mathcal B(x,r)\cap \Omega,
|y|^\lambda \ d^\alpha (y) dy )) \cap C^0((0,r^2); H^1(\mathcal
B(x,r)\cap \Omega, |y|^\lambda \ d^\alpha (y) dy))$ such that for
any $\Phi\in C^0((0,r^2);C^\infty_0(\mathcal B(x,r) \cap \Omega))$
and any $0<t_1<t_2<r^2$ we have \be \int_{t_1}^{t_2}\int_{\mathcal
B(x,r)\cap \Omega} \{ |y|^\lambda d^\alpha (y)v_t \Phi+|y|^\lambda
d^\alpha(y) \nabla v\nabla \Phi\} dy dt = 0 \ .\ee
\end{definition}

Then we have

\begin{theorem}
\label{harnackgen} Let $\alpha\ge 1$, $N\ge 2$,
$\lambda\in[2-N,0]$ and $\Omega\subset \R^N$ be a smooth bounded
domain containing the origin. Then there exist positive constants
$C_H$ and $R=R(\Omega)$ such that for $x\in \Omega$, $0<r<R$ and
for any positive solution $v(y,t)$ of $\frac{\partial v}{\partial
t}=\frac{1}{|y|^\lambda d^\alpha(y)} div (|y|^\lambda d^\alpha (y)
\nabla v)$ in $\left\{\mathcal B(x,r)\cap \Omega\right\}\times
(0,r^2)$, the following estimate holds true
$${\it ess~sup}_{(y,t)\in \left\{\mathcal B(x,\frac{r}{2})\cap \Omega
\right\}\times (\frac{r^2}{4},\frac{r^2}{2})} v(y,t) \le C_H ~
{\it ess~inf}_{(y,t)\in \left\{\mathcal B(x,\frac{r}{2})\cap
\Omega \right\}\times (\frac{3}{4} r^2,r^2)} v(y,t)\ .$$
\end{theorem}

\smallskip

In order to prove the parabolic Harnack inequality in Theorem
\ref{harnackgen} we use the Moser iteration technique as adapted
to degenerate elliptic operators in \cite{FKS}, \cite{CS} as well
as \cite{GSC}. In this approach one inserts in the weak form of
the equation $v_t=-L^\lambda_\alpha v$ suitable test functions
$\Phi$. One of the key ideas is to use test functions $\Phi$ of
the form $\eta^2 v^q$, where $v$ is the weak solution of the
equation, $\eta$ is a cut off function and $q \in \R$. To this end
one has to check that $\eta^2 v^q$ is in the right space of test
function. In this direction the following density theorem is
crucial.

\begin{theorem}
\label{density} Let $N\ge 2$ and $\Omega\subset \R^N$ be a smooth
bounded domain. Then for any $\alpha\ge 1$
$$H^1(\Omega, \ d^\alpha(y) \ dy)=H^1_0(\Omega, d^\alpha(y) \ dy) \ .$$
In particular for any $\alpha \ge 1$, the set $C^\infty_0(\Omega)$
is dense in $H^1(\Omega, d^\alpha(y) dy)$.
\end{theorem}

Here $H^{1}(\Omega, d^\alpha(y)dy)$ denotes the set $\{v=v(y):
\int_{\Omega} d^\alpha (y) (v^2+|\nabla v|^2)  dy<\infty\}$, the
corresponding norm being defined in (\ref{ehi}).

\medskip

We are now ready to prove the density theorem.

\smallskip
{\em Proof } Let us prove here the result when $\alpha=1$. We
refer to Proposition 9.10 in \cite{K} for the case $\alpha>1$,
even though our proof with some minor changes, can also cover this
range.

First of all from Theorem 7.2 in \cite{K} it is known that the set
$C^\infty(\overline \Omega)$ is dense in $H^1(\Omega, d(y) \ dy)$.
Thus for any $v\in H^1(\Omega, d(y) \ dy)$ there exists $v_m\in
C^\infty(\overline \Omega)$ such that for any $\epsilon>0$ we have
$||v-v_m||_{H^1_{1}}\le \epsilon$ if $m\ge m(\epsilon)$. Let us
choose $w:=v_{m(\epsilon)}$ and let us define, for $k\ge 1$, the
following function
$$\varphi_k(x)=\begin{cases} 0 & if\  d(x)\le \frac{1}{k^2}  ~ ,\\
1 +\frac{\ln (kd(x))}{\ln (k)} & if \ \frac{1}{k^2}< d(x)<
\frac{1}{k}~ , \\
 1 & if \ d(x)\ge \frac{1}{k} \ . \end{cases} $$
Then $w_k:=w \varphi_k\in C^{0,1}_0(\Omega)$, moreover we have
$$||w-w_k||_{H^1_1}= ||w
(1-\varphi_k)||_{H^1_1}\le 2\int_{\Omega} (w^2+|\nabla w|^2)
(1-\varphi_k)^2 d(y) \ dy + 2 \int_{\Omega} w^2 |\nabla
\varphi_k|^2 d(y) \ dy \le$$ $$\le  2\int_{d(y)< \frac{1}{k}}
(w^2+|\nabla w|^2) d(y) \ dy + 2 \int_{\frac{1}{k^2}<d(y)<
\frac{1}{k}} \frac{w^2}{d(y) (\ln (k))^2} \ dy \ .$$ Now as $k\to
\infty$ the right hand side goes to zero, this proves the Theorem.

\finedim

The above Theorem allows us to take the cut off function $\eta$ in
$C^\infty_0(\mathcal B(x,r))$ instead of taking it as usual in
$C^\infty_0(\mathcal B(x,r)\cap \Omega)$. Clearly the two function
spaces differ only if the ``ball" intersects the boundary of
$\Omega$. To explain what are the appropriate modifications of the
standard iteration argument by Moser, we now present in detail the
first step, which is the $L^2$ mean value inequality for any
positive local subsolution of the equation $v_t=-L^\lambda_\alpha
v$.

\begin{theorem}
\label{meanval} Let $\alpha\ge 1$, $N\ge 2$, $\lambda \in[2-N,0]$
and $\Omega\subset \R^N$ be a smooth bounded domain containing the
origin. Then there exist positive constants $C$ and $R(\Omega)$
such that for $x\in \Omega$, $0<r<R(\Omega)$ and for any positive
subsolution $v(y,t)$ of $v_t-\frac{1}{|y|^\lambda d^\alpha(y)}
div(|y|^\lambda d^\alpha(y) \nabla v)=0$ in $\{\mathcal B(x,r)\cap
\Omega\}\times (0,r^2)$ we have the estimate
$${\it ess~ sup}_{(y,t)\in \{\mathcal B(x,\frac{r}{2})\cap \Omega\}\times (\frac{r^2}{2}, r^2)} \ \ v^2(y,t)\le
\frac{C}{r^2 V(x,r)} \int_{\{\mathcal B(x,r)\cap \Omega\}\times
(0,r^2)} |y|^\lambda d^\alpha(y) \  v^2(y,t)dy dt \ .$$
\end{theorem}

{\em Proof: } We will only prove the result in the non standard
case in which the ``ball" $\mathcal B(x,r)$ intersects the
boundary of $\Omega$; we refer to \cite{MT2} as well as to
\cite{GSC} for details in the other case. Similarly to Definition
\ref{sol} we define a subsolution $v(y,t)$ to be a function in
$C^1((0,r^2);L^2(\mathcal B(x,r)\cap \Omega, |y|^\lambda
d^\alpha(y) \ dy ))\cap C^0((0,r^2);H^1(\mathcal B(x,r)\cap
\Omega, |y|^\lambda d^\alpha (y) \ dy))$ such that the following
holds true
$$\int^{r^2}_{0}\int_{\mathcal B(x,r)\cap \Omega} \{|y|^\lambda d^\alpha (y) v_t\Phi+|y|^\lambda
d^\alpha(y) \nabla v\nabla
\Phi\} dy dt \le 0, ~ ~ ~ \forall \ \Phi \in
C^0((0,r^2);C^\infty_0(\mathcal B(x,r)\cap \Omega)) \ , \Phi\ge
0.$$ Hence in particular we have also
$$\int_{\mathcal B(x,r)\cap \Omega} \{|y|^\lambda d^\alpha (y) v_t\Phi+|y|^\lambda
d^\alpha(y) \nabla v\nabla \Phi\} dy \le 0, ~ ~ ~ \forall \ \Phi
\in C^\infty_0(\mathcal B(x,r)\cap \Omega) \ , \Phi\ge 0.$$ Let us
define for any $q, M \ge 1$ the following functions $G(z)=z^q$ if
$z\le M$ and $G(z)=M^q+q(z-M)M^{q-1}$ if $z> M$ and $H(z)\ge 0$ by
$H'(z)=\sqrt{G'(z)}$, $H(0)=0$; note that $G(z)\le zG'(z)$ as well
as $H(z)\le zH'(z)$. Due to Theorem \ref{density} there exists a
sequence of functions $v_m$ in $C^\infty(\overline{\mathcal
B(x,r)\cap \Omega})$ having compact support in $\Omega$ such that
$v_m\to  v$ in $H^1(\mathcal B(x,r)\cap \Omega, d^\alpha (y) \
dy)$ as $m\to +\infty$; whence due to (\ref{ohi}) also in
$H^1(\mathcal B(x,r)\cap \Omega, |y|^\lambda d^\alpha (y) \ dy)$.
Hence for any $\eta\in C^\infty_0(\mathcal B(x,r))$ and $m\ge 1$
the function $\Phi:=\eta^2 G(v_m)$ is an admissible test function,
that is the following holds true
$$\int_{\mathcal B(x,r)\cap \Omega} \{|y|^\lambda d^\alpha(y) \eta^2 G(v_m) v_t +|y|^\lambda d^\alpha(y) \nabla v\nabla
(\eta^2 G(v_m))\} dy \le 0 \ .$$ Passing to the limit as $m\to
+\infty$ we get
$$\int_{\mathcal B(x,r)\cap \Omega} \{ |y|^\lambda d^\alpha(y) \eta^2 G(v) v_t+ |y|^\lambda
 d^\alpha(y) \nabla v\nabla
(\eta^2 G(v))\} dy\le 0 \ , ~ ~ ~\forall \ \eta\in
C^\infty_0(\mathcal B(x,r))\ .$$ This is the standard starting
point in Moser iteration technique apart from the fact that the
cut off function $\eta$ does not be necessarily zero on $\partial
\Omega$, this is crucial. Then by Schwarz inequality we get
$$\int_{\mathcal B(x,r)\cap \Omega} \{ |y|^\lambda d^\alpha(y) \eta^2 G(v) v_t+ |y|^\lambda d^\alpha(y) |\nabla v|^2
 G'(v) \eta^2 \} dy \le
C \int_{\mathcal B(x,r)\cap \Omega} |y|^\lambda d^\alpha(y)
|\nabla \eta|^2 v^2 G'(v) dy $$ thus also that
$$\int_{\mathcal B(x,r)\cap \Omega} \{|y|^\lambda d^\alpha(y) \eta^2 G(v) v_t+ |y|^\lambda d^\alpha(y)
|\nabla (\eta H(v)) |^2 \} dy \le C \int_{\mathcal B(x,r)\cap
\Omega} |y|^\lambda d^\alpha(y) |\nabla \eta|^2 v^2 G'(v) dy \ .$$
For any smooth function $\chi$ of the time variable $t$, we easily
get $$\frac{d}{dt}\int_{\mathcal B(x,r)\cap \Omega}|y|^\lambda
d^\alpha(y) (\eta \chi F(v))^2 dy+ \chi^2 \int_{\mathcal
B(x,r)\cap \Omega} |y|^\lambda d^\alpha(y) |\nabla(\eta H(v))|^{2}
dy \le$$
$$\le C \chi \left(\chi ||\nabla \eta||_{L^\infty(\R^n)}
+||\chi'||_{L^\infty(\R)}\right) \int_{{\it supp} \eta\ \cap
\Omega} |y|^\lambda d^\alpha(y) v^2 G'(v) dy ;
 $$ here $F(z)$ is such that $2F(z)F'(z)=G(z)$.
For $\frac{1}{2}\le s<s'<1$ we choose as usual $\chi$ such that
$0\le \chi\le 1$, $\chi=0$ in $(-\infty, r^2(1-s'))$, $\chi=1$ in
$(r^2(1-s), \infty)$, moreover if $\xi \in C^\infty_0(0,1)$ be a
nonnegative non increasing function such that $\xi(z)=1$ if $z\le
s$ and $\xi(z)=0$ if $z\ge s'$, we define, making use of local
coordinates, the following cut off function
$\eta(y):=\xi\left(\frac{|y'-x'|}{r}\right)\xi\left(\frac{|a(y')-y_N-d(x)|}{r}\right)$.
Then clearly $||\nabla \eta||_{L^\infty(\R^n)}\le
\frac{C}{r(s'-s)}$ and $||\chi'||_{L^\infty(\R)}\le
\frac{C}{r^2(s'-s)}$.

Integrating our inequality over $(0,t)$, with $t\in(r^2(1-s),
r^2)$ we obtain
$$\sup_{t\in J}\int_{\mathcal B(x,r)\cap \Omega} |y|^\lambda d^\alpha(y)
(\eta F(v))^2 dy +\int_{\{\mathcal B(x,r)\cap \Omega\}\times
(r^2(1-s), r^2)} |y|^\lambda d^\alpha(y) |\nabla(\eta H(v))|^{2}
dy dt \le$$
$$\le
\frac{C}{r^2(s'-s)^2}\int_{\{\mathcal B(x,s'r) \cap \Omega\}\times
(r^2(1-s'),r^2)} |y|^\lambda d^\alpha(y) v^2 G'(v) dy dt \ .$$
Making once again use of Theorem \ref{density} we note that we can
apply the local weighted Moser inequality in Theorem \ref{moser}
to the function $f:=\eta H(v)$ thus obtaining
$$\int_{\{\mathcal B(x,r)\cap \Omega\}\times (r^2(1-s),
r^2)} |y|^\lambda d^\alpha(y) (\eta
H(v))^{2\left(1+\frac{2}{N+\alpha}\right)} dy dt\le $$
$$\le\frac{C}{r^2(s'-s)^2} \left(\int_{\{\mathcal B(x,s'r)\cap
\Omega\}\times (r^2(1-s'),r^2)} |y|^\lambda d^\alpha(y) v^2 G'(v)
dy dt \right)^{1+\frac{2}{N+\alpha}} \ .$$ Let us now denote by
$\tilde\gamma:=1+\frac{2}{N+\alpha}$ thus as $M$ tends to infinity
we have for $p:=q+1$
$$\int_{\{\mathcal B(x,sr)\cap \Omega\}\times (r^2(1-s),r^2)}|y|^\lambda  d^\alpha(y) v^{p\tilde\gamma} dy dt\le
\frac{C}{r^2(s'-s)^2} \left(\int_{\{\mathcal B(x,s'r)\cap
\Omega\}\times (r^2(1-s'),r^2)} |y|^\lambda d^\alpha(y) v^p dy dt
\right)^{\tilde\gamma} \ .$$ Thus due to Lemma \ref{volume} also
that
$$V(x,sr)^{-1} (r^2s)^{-1}\int_{\{\mathcal B(x,sr)\cap \Omega\}\times(r^2(1-s),r^2)} |y|^\lambda d^\alpha(y)
v^{p\tilde \gamma} dy dt\le$$ $$\le C
\left(\frac{s}{s'-s}\right)^2 \left(V(x,s'r)^{-1}
(r^2s')^{-1}\int_{\{\mathcal B(x,s'r)\cap
\Omega\}\times(r^2(1-s'),r^2)} |y|^\lambda d^\alpha(y) v^p dy
dt\right)^{\tilde \gamma} \ .$$ Take now $p=p_i:=2{\tilde
\gamma}^{i}$, $s=\theta_{i+1}$ and $s'=\theta_i$ where
$\theta_i:=\frac{i+2}{2(i+1)}$ then if we denote by
$I(i):=\left(V(x,\theta_i r)^{-1} (r^2\theta_i)^{-1}
\int_{\{\mathcal B(x,\theta_i r)\cap \Omega\}\times
(r^2(1-\theta_i),r^2)} |y|^\lambda d^\alpha(y) v^{p_i} dy
dt\right)^{\frac{1}{p_i}}$ the above inequality can be restated as
follows $I(i+1)\le C(i) I(i)$. Thus since one can show that the
product of $C(i)$ for all $i\ge 0$ is finite, we obtain
$I(\infty)\le \left\{\prod^\infty_{i=0} C(i)\right\} I(0)$, this
completes the proof of the proposition. To this end the choice
$R(\Omega):=\min\{\beta,R(1,\Omega)\}$ can be made, here $\beta$
and $R(1,\Omega)$ are the constants appearing respectively in the
local representation of $\partial \Omega$ and in Theorem
\ref{moser} when $\alpha:=1$.

\finedim

Theorem \ref{moser} corresponds to the local weighted Moser
inequality needed in the proof of the parabolic Harnack inequality
up to the boundary stated in Theorem \ref{harnack}. The local
weighted Moser inequality involved in the proof of Theorem
\ref{harnackgen} differs from Theorem \ref{moser} only if $d(x)\ge
\gamma r$, $N\ge 3$, $\lambda \neq 0$, and in this case it reads
as follows

\begin{theorem}
\label{ultima} Let  $N\ge 3$, $\lambda\in [2-N,0)$ and
$\Omega\subset \R^N$ be a smooth bounded domain containing the
origin. Then there exist a positive constant $C_M$ such that for
any $\nu\ge N$, $x\in \Omega$, $r>0$ and $f \in
C^\infty_0(\mathcal B(x,r))$ we have
$$\int_{B(x,r)}|y|^\lambda |f(y)|^{2\left(1+\frac{2}{\nu}\right)}  dy\le C_M r^2
(r^N (|x|+r)^\lambda)^{-\frac{2}{\nu}}
\left(\int_{B(x,r)}|y|^\lambda |\nabla f|^2 dy \right)
\left(\int_{B(x,r)} |y|^\lambda |f|^2 dy \right)^\frac{2}{\nu} .
$$
\end{theorem}

{\em Proof :} By H\"older inequality the result easily follows
with $C_M:=C_S$ as soon as the following local weighted Sobolev
inequality holds true \be\label{ls2}\left(\int_{B(x,r)}|y|^\lambda
|f(y)|^{\frac{2N}{N-2}}  dy\right)^{\frac{N-2}{N}}\le C_S
(|x|+r)^{\frac{2|\lambda|}{N}}\int_{B(x,r)}|y|^\lambda |\nabla
f|^2 dy \ee (we refer to the proof of Theorem \ref{moser} where a
similar argument is used). Let us first prove the above inequality
for any $\lambda\in (2-N,0)$. As a consequence of the Caffarelli
Kohn Nirenberg inequality (e.g. see Corollary 2 in Section 2.1.6
of \cite{M}), the following holds true
$$\left(\int_{B(x,r)} f^{\frac{2N}{N-2}}
|y|^{\frac{N\lambda}{N-2}} dy \right)^\frac{N-2}{N} \le C
\int_{B(x,r)} |\nabla f|^2 |y|^\lambda dy \ , ~ ~ ~ \forall f\in
C^\infty_0(B(x,r))$$ and for some positive constant $C$
independent of $x$ and $r$. Whence also that
$$\left(\int_{B(x,r)} f^{\frac{2N}{N-2}}
|y|^{\lambda} dy \right)^{\frac{N-2}{N}} \le C \left(\sup_{y\in
B(x,r)} |y|\right)^{\frac{2|\lambda|}{N}} \int_{B(x,r)} |\nabla
f|^2 |y|^\lambda dy \le C (|x|+r)^{\frac{2|\lambda|}{N}}
\int_{B(x,r)} |\nabla f|^2 |y|^\lambda dy \ .
$$
Let us now prove the result for $\lambda=2-N$. To this end let us
apply Proposition \ref{prop6} to $\Omega=B(0,1)$ with
$D=e^{\frac{1}{N-2}}$. Then there exists a positive constant $C$
such that
$$\int_{B(0,1)} |\nabla v|^2 |x|^{2-N} dx \ge C \left(\int_{B(0,1)} v^{\frac{2N}{N-2}}
|x|^{-N} X^{\frac{2(N-1)}{N-2}}\left(\frac{|x|}{D}\right) dx
\right)^{\frac{N-2}{N}}\ , ~ ~ ~ \forall \ v\in
C^\infty_0(B(0,1))\ ;$$ here $X(t)=\frac{1}{1-\ln t}$, $t\in
(0,1]$. Now let us take $v(x):=f\left(\frac{x}{R}\right)$ for any
$f\in C^\infty_0(B(0,R))$ then from above we have
$$\int_{B(0,R)} |\nabla f|^2 |y|^{2-N} dy \ge C \left(\int_{B(0,R)} f^{\frac{2N}{N-2}}
|y|^{-N} X^{\frac{2(N-1)}{N-2}}\left(\frac{|y|}{DR}\right) dy
\right)^{\frac{N-2}{N}}\ .$$ Then if $y\in B(x,r)$ clearly $y\in
B(0,|x|+r)$, thus if we take $R:=|x|+r$ and $f\in
C^\infty_0(B(x,r))$ from above we have
$$\int_{B(x,r)} |\nabla f|^2 |y|^{2-N} dy \ge C \left(\int_{B(x,r)} f^{\frac{2N}{N-2}}
|y|^{-N} X^{\frac{2(N-1)}{N-2}}\left(\frac{|y|}{DR}\right) dy
\right)^{\frac{N-2}{N}}\ge$$ $$\ge  \left(\int_{B(x,r)}
f^{\frac{2N}{N-2}} |y|^{2-N} dy \right)^{\frac{N-2}{N}} \left
(\inf_{y\in B(x,r)} |y|^{-2}
X^{\frac{2(N-1)}{N-2}}\left(\frac{|y|}{DR}\right)
\right)^{\frac{N-2}{N}} \ .$$ Whence the claim easily follows as
soon as we prove that $$\left(\sup_{y\in B(x,r)} |y|
X^{-\frac{N-1}{N-2}}\left(\frac{|y|}{DR}\right)
\right)^{\frac{2(N-2)}{N}}\le C_S (|x|+r)^{\frac{2(N-2)}{N}}\ .$$
This is indeed the case in fact we have
$$\sup_{y\in B(x,r)} |y|
X^{-\frac{N-1}{N-2}}\left(\frac{|y|}{DR}\right) \le
\sup_{0\le|y|\le|x|+r} |y| \left(1-\ln\left(\frac{|y|}{DR}\right)
 \right)^{\frac{N-1}{N-2}} =$$
(thus using the fact that the function
$\varphi(t)=t\left(1-\ln\left(
\frac{t}{DR}\right)\right)^{\frac{N-1}{N-2}}$ is an increasing
function for $t\in [0,R]$ if $D$ and $R$ are as above)
$$= (|x|+r) \left(1-\ln\left(\frac{|x|+r}{DR}\right)
 \right)^{\frac{N-1}{N-2}} =(|x|+r) (1+\ln(D))^{\frac{N-1}{N-2}}=(|x|+r)
 \left(\frac{N-1}{N-2}\right)^{\frac{N-1}{N-2}} \ .
$$
This completes the proof of Theorem \ref{ultima}.

\finedim

\medskip

To state the heat kernel estimates following from Theorem
\ref{harnackgen} we introduce some notation. The operator
$L^\lambda_\alpha$ is defined for $\alpha \ge 1$ and $\lambda\in
[2-N,0]$ in $L^2(\Omega, |x|^\lambda d^\alpha(x) \ dx )$ as the
generator of the symmetric form
$$\mathcal L^\lambda_\alpha[v_1,v_2]:=\int_{\Omega} |x|^\lambda d^\alpha(x) \nabla v_1
 \nabla v_2 \ dx \ ,$$
namely
$$D(L^\lambda_\alpha):=\left\{v\in H^1_0(\Omega, |x|^\lambda d^\alpha(x) \ dx ):
-\frac{1}{|x|^\lambda
 d^\alpha(x)} div (|x|^\lambda d^\alpha(x) \nabla v)
\in L^2(\Omega, |x|^\lambda d^\alpha(x) \ dx)\right\},$$
$$L^\lambda_\alpha v:=  -\frac{1}{|x|^\lambda d^\alpha(x)} div (|x|^\lambda d^\alpha(x) \nabla v)
\hbox { for any } v\in D(L^\lambda_\alpha) \ ,$$ where
$H^1_0(\Omega,|x|^\lambda d^\alpha(x) \ dx )$ denotes the closure
of $C^\infty_0(\Omega)$ in the norm
\begin{equation}
\label{ehibis}
 v\to ||v||_{H^1_{\alpha, \lambda}}:=\left\{\int_{\Omega}
|x|^\lambda d^\alpha(x) \left(|\nabla v|^2 +v^2\right) \ \
dx\right\}^{\frac{1}{2}} \ .
\end{equation}

Then $L^\lambda_\alpha $ is a nonnegative self-adjoint operator on
$L^2(\Omega, |y|^\lambda d^\alpha(y) dy)$ such that for every
$t>0$, $e^{-L^\lambda_\alpha t}$ has a integral kernel, that is
$e^{-L^\lambda _\alpha t}v_0(x):=\int_{\Omega} l^\lambda_\alpha
(t,x,y)v_0(y) |y|^\lambda d^\alpha(y) dy$; here $l^\lambda_\alpha
(t,x,y)$ is called the heat kernel of $L^\lambda_\alpha$. The
existence of $l^\lambda_\alpha (t,x,y)$ can be proved arguing as
in \cite{DS1}; that is, using a global Sobolev inequality on
$\Omega$, which can be easily deduced from its local version
(\ref{ls}) as well as (\ref{ls2}), by means of the partition of
unity as in \cite{K}.

Then from the parabolic Harnack inequality in Theorem
\ref{harnackgen} the following sharp two-sided heat kernel
estimate can be easily deduced:

\begin{theorem}
\label{heatgen} Let $\alpha\ge 1$, $N\ge 2$, $\lambda\in [2-N,0]$
and $\Omega\subset \R^N$ be a smooth bounded domain containing the
origin. Then there exist positive constants $C_1, C_2$, with
$C_1\le C_2$, and $T>0$ depending on $\Omega$ such that
$$C_1 \min\left\{\frac{1}{t^{\frac{\alpha}{2}}},
\frac{(|x|+\sqrt t)^{\frac{|\lambda|}{2}}(|y|+\sqrt
t)^{\frac{|\lambda|}{2}}}{d^{\frac{\alpha}{2}}(x)d^{\frac{\alpha}{2}}(y)}\right\}
t^{-\frac{N}{2}} e^{-C_2\frac{|x-y|^2}{t}}\le
l^\lambda_\alpha(t,x,y) \le $$ $$\le C_2
\min\left\{\frac{1}{t^{\frac{\alpha}{2}}},\frac{(|x|+\sqrt
t)^{\frac{|\lambda|}{2}}(|y|+\sqrt
t)^{\frac{|\lambda|}{2}}}{d^{\frac{\alpha}{2}}(x)
d^{\frac{\alpha}{2}}(y)}\right\}t^{-\frac{N}{2}}
e^{-C_1\frac{|x-y|^2}{t}}$$ for all $x,y\in \Omega$ and $0<t\le
T$.
\end{theorem}

\smallskip
{\em Proof of Theorem \ref{heatgen}: } Using the mean value
estimate for subsolutions as in Theorem \ref{meanval} and the
parabolic Harnack inequality of Theorem \ref{harnackgen} and
arguing as in Theorems 5.2.10, 5.4.10 and 5.4.11 in \cite{SC2} we
are lead to the following Li-Yau type estimate
$$\frac{C_1~  e^{-C_2\frac{|x-y|^2}{t}}}{V(x,\sqrt t)^{\frac{1}{2}}
V(y,\sqrt t)^{\frac{1}{2}}} \le l^\lambda_\alpha (t,x,y) \le
\frac{C_2~ e^{-C_1\frac{|x-y|^2}{t}}}{V(x,\sqrt t)^{\frac{1}{2}}
V(y,\sqrt t)^{\frac{1}{2}}} \ ,$$ for all $x,y\in \Omega$ and
$0<t\le T$; where $C_1, C_2$ are two positive constants with
$C_1\le C_2$, and $T>0$ depends on $\Omega$. From this the result
follows using the volume estimate in Lemma \ref{volume}.

\finedim

Using all the machinery we have produced in this section we can
handle more general operators than the one in Theorems
\ref{harnackgen} and \ref{heatgen}. Thus, consider the operator
\be \widetilde{L^\lambda_\alpha} := -\frac{1}{|x|^\lambda
d^\alpha(x)}\sum^N_{i,j=1} \frac{\partial }{\partial x_i}\left(
a_{i,j}(x) |x|^\lambda d^\alpha(x) \frac{\partial}{\partial
x_j}\right) \ , \ee where $\left(a_{i,j}(x)\right)_{n\times n}$ is
a measurable symmetric uniformly elliptic matrix. The operator
$\widetilde{L^\lambda_\alpha}$ is defined for $\alpha \ge 1$ and
$\lambda\in [2-N,0]$ in $L^2(\Omega, |x|^\lambda d^\alpha(x) \ dx
)$ as the generator of the symmetric form
$$\mathcal {\widetilde{L^\lambda_\alpha}}[v_1,v_2]:=\sum^N_{i,j=1} \int_{\Omega}
|x|^\lambda d^\alpha(x) a_{i,j}(x)\frac{\partial v}{\partial
x_i}\frac{\partial v}{\partial x_j} dx \ .$$ Then the existence of
a heat kernel $\widetilde{l^\lambda_\alpha}(t,x,y)$ follows as in
\cite{DS1}, and we have

\begin{theorem}
\label{stab} Let $\alpha\ge 1$, $N\ge 2$, $\lambda\in [2-N,0]$ and
$\Omega\subset \R^N$ be a smooth bounded domain containing the
origin. Then there exist positive constants $C_1, C_2$, with
$C_1\le C_2$, and $T>0$ depending on $\Omega$ such that
$$C_1 \min\left\{\frac{1}{t^{\frac{\alpha}{2}}},
\frac{(|x|+\sqrt t)^{\frac{|\lambda|}{2}}(|y|+\sqrt
t)^{\frac{|\lambda|}{2}}}{d^{\frac{\alpha}{2}}(x)d^{\frac{\alpha}{2}}(y)}\right\}
t^{-\frac{N}{2}} e^{-C_2\frac{|x-y|^2}{t}}\le
\widetilde{l^\lambda_\alpha}(t,x,y) \le $$ $$\le C_2
\min\left\{\frac{1}{t^{\frac{\alpha}{2}}},\frac{(|x|+\sqrt
t)^{\frac{|\lambda|}{2}}(|y|+\sqrt
t)^{\frac{|\lambda|}{2}}}{d^{\frac{\alpha}{2}}(x)
d^{\frac{\alpha}{2}}(y)}\right\}t^{-\frac{N}{2}}
e^{-C_1\frac{|x-y|^2}{t}}$$ for all $x,y\in \Omega$ and $0<t\le
T$.
\end{theorem}

\begin{remark}\label{stab2} A parabolic Harnack inequality up to the boundary similar to the one of
Theorem \ref{harnackgen} can be stated under the same assumptions
of Theorem \ref{stab} for the more general operator $\widetilde
{L^\lambda_\alpha}$.
\end{remark}

\setcounter{equation}{0}
\section{Critical point singularity}

In this section we  establish a new Improved Hardy inequality
(Theorem \ref{prop7}) and then we give the proofs of Theorem
\ref{thmorigin} and Theorem \ref{thmoriginbis}. The structure of
this section is as follows.

In Subsection 3.1 we first deduce  the improved Hardy inequality
and then the global in time pointwise upper bound for the heal
kernel of the Schr\"odinger operator $-\Delta
-((N-2)^2/4)|x|^{-2}$, which is sharp when $x$ and $y$ are close
to the boundary (see Theorem \ref{th5}); then, due to an argument
contained in \cite{D1}, we complete the proof of Theorem
\ref{thmoriginbis} proving the sharp lower bound for time large
enough.

The proof of Theorem \ref{thmorigin} is finally completed in
Subsection 3.2, using the parabolic Harnack inequality up to the
boundary of Theorem \ref{harnackgen}.

\subsection {Boundary upper bounds and complete sharp
 description of the heat kernel for large values of time}

We first recall the following improved Hardy-Sobolev inequality
stated in Theorem A in \cite{FT} (see also inequality (3.3) in
\cite{BFT2})

\begin{proposition}
\label{prop6} For $N\ge 3$, let $\Omega\subset \R^N$ be a smooth
bounded domain containing the origin and $D\ge \sup_{x\in \Omega}
|x|$. Then there exists a positive constant $C$ such that
$$\int_\Omega |\nabla v|^2 |x|^{2-N} dx \ge C \left(\int_{\Omega} v^{\frac{2N}{N-2}}
|x|^{-N} X^{\frac{2(N-1)}{N-2}}\left(\frac{|x|}{D}\right) dx
\right)^{\frac{N-2}{N}}\ ,$$ for any $v\in C^\infty_0(\Omega)$;
here $X(t)=\frac{1}{1-\ln t}$, $t\in (0,1]$.
\end{proposition}
We next state a new result, the proof of which will be given later on.
\begin{theorem}
{\bf (Improved Hardy inequality)} \label{prop7} Let $\Omega\subset
\R^N$, $N\ge 3$, be a smooth bounded domain containing the origin.
Then there exists a constant $C=C(\Omega)\in (0,\frac{1}{4}]$ such
that
\begin{equation} \label{improvedHardy}\int_\Omega \left(|\nabla u|^2 -\frac{(N-2)^2}{4 |x|^2} u^2 \right) dx
 \ge C(\Omega) \int_{\Omega} \frac{u^2}{d^2(x)} \  dx  \ , ~ ~ ~ \forall \ u \in C^\infty_0(\Omega)\ . \end{equation}
The positive constant $C(\Omega)$ can be taken to be exactly
$\frac{1}{4}$ for all domains satisfying the following condition
\begin{equation}
\label{cond} -div (|x|^{2-N} \nabla d(x))\ge 0 \ \hbox { a.e. in }
\Omega \ .
\end{equation}
\end{theorem}
For example when $\Omega\equiv B(0,R)$, for  arbitrary $R>0$,
condition (\ref{cond}) is satisfied. Consequently, in this case
the Hardy inequality involving the Schr\"odinger operator having
critical singularity at the origin can be improved exactly by the
inverse-square potential having critical singularity at the
boundary.

As a consequence of Proposition \ref{prop6} and of the improved
Hardy inequality of Theorem \ref{prop7}, the following logarithmic
Sobolev inequality can be easily obtained:

\begin{theorem} {\bf (Logarithmic Hardy Sobolev inequality)}
\label{teorema} For $N\ge 3$, let $\Omega\subset \R^N$ be a smooth
bounded domain containing the origin. Then for any $u\in
C^\infty_0(\Omega\setminus \{0\})$, $u\ge 0$, and any $\epsilon
>0$ we have
\begin{equation} \label{logsob3} \int_{\Omega}
u^2 \log \left(\frac{u}{||u||_2 |x|^{\frac{2-N}{2}} d(x)}\right)
dx \le \epsilon \int_{\Omega}\left(|\nabla u|^2
-\frac{(N-2)^2}{4|x|^2}  u^2 \right)dx  +\left(K_3 -\frac{N+2}{4}
\log \epsilon \right)||u||^2_2 \ ;
\end{equation}
here $K_3$ is a positive constant independent of $\epsilon$ and
$||u||_2:=\left(\int_{\Omega} u^2 dx\right)^{\frac{1}{2}}$.
\end{theorem}
\medskip

Then using Gross theorem of logarithmic Sobolev inequalities, as
adapted by Davies and Simon (see Theorem 2.2.7 in \cite{D4}), we
will show the following global in time pointwise upper bound for
the heat kernel:

\begin{theorem}
\label{th5} For $N\ge 3$, let $\Omega \subset \R^N$ be a smooth
bounded domain containing the origin. Then there exists a positive
constant $C$ such that
$$k(t,x,y)\le C \frac{d(x) d(y)}{t} |xy|^{\frac{2-N}{2}}  t^{-\frac{N}{2}} e^{-\lambda_1 t}
 , ~ ~ ~ \forall \ x,y\in \Omega , \ t>0 \ .$$
\end{theorem}

Let us first prove the logarithmic Hardy Sobolev inequality
(\ref{logsob3}).

\medskip
{\em Proof of Theorem \ref{teorema}: } As a first step we claim
that the following logarithmic Hardy Sobolev inequality holds
true:

\begin{equation}
\label{logsob1}\int_{\Omega} u^2 (-\log d(x)) \ dx \le \epsilon \int_{\Omega}
\left(|\nabla u|^2 -\frac{(N-2)^2}{4|x|^2}  u^2 \right) dx  +\left(K_1-\frac{1}{2} \log \epsilon\right) ||u||^2_2\ ,
\end{equation}
for any $u\in C^\infty_0(\Omega)$, $u\ge 0$, and any $\epsilon>0$;
here $K_1$ is a positive constant independent of $\epsilon$.

To see this let us first suppose that the nonnegative function
$u\in C^{\infty}_0(\Omega)$ is such that $||u||_2=1$. We then have

$$\int_{\Omega} u^2 (-\log d(x))\ dx = \frac{1}{2} \int_{\Omega} u^2 (\log d(x)^{-2})\ dx
 \le \frac{1}{2} \log \left(\int_{\Omega} \frac{1}{d(x)^2} u^2 dx\right) \le $$
$$\le \frac{1}{2}\log \left(C^{-1}\int_{\Omega} \left(|\nabla u|^2 -\frac{(N-2)^2}{4|x|^2} u^2
\right) dx \right) \ ;$$ here we have used first Jensen's
inequality and then the improved Hardy inequality
(\ref{improvedHardy}). For a general nonnegative $u\in
C^\infty_0(\Omega)$ we apply the above inequality to the function
$\frac{u}{||u||_2}$, to get
$$\int_{\Omega} u^2 (-\log d(x))\ dx \le
\frac{1}{2}||u||^2_2 \log
\left(\frac{C^{-1}}{||u||^2_2}\int_{\Omega} \left(|\nabla u|^2
-\frac{(N-2)^2}{4|x|^2} u^2 \right) dx \right) \ .$$ Since $\log
z\le z$ for any $z>0$, then also $\log y\le \epsilon 2 C y -\log
\left(\epsilon 2 C\right)$, for any $\epsilon>0$; whence from this
we deduce (\ref{logsob1}), with $K_1:= \frac{1}{2}\log
(\frac{1}{2C})$.

We will next show the following logarithmic Hardy Sobolev
inequality:

\begin{equation}
\label{logsob2}
\int_{\Omega} u^2 \log \left(\frac{u}{||u||_2 |x|^{\frac{2-N}{2}}}\right) dx \le
\epsilon \int_{\Omega}\left(|\nabla u|^2 -\frac{(N-2)^2}{4|x|^2}  u^2 \right)dx  +\left(K_2 -\frac{N}{4}
 \log \epsilon \right)||u||^2_2 \ ,
\end{equation}
for any $u\in C^\infty_0(\Omega \setminus \{0\})$, $u\ge 0$, and
any $\epsilon >0$; here $K_2$ is a positive constant independent
of $\epsilon$.

By Proposition \ref{prop6} it follows easily  that there exists a
positive constant $C$ such that

\begin{equation}
\label{ft} \int_\Omega |\nabla v|^2 |x|^{2-N} dx \ge C
\left(\int_{\Omega} v^{\frac{2N}{N-2}} |x|^{2-N} dx
\right)^{\frac{N-2}{N}} ,
\end{equation}
for any $v\in C^\infty_0(\Omega)$ (this is inequality (4.12) in
\cite{BFT2}). Whence we claim that the following logarithmic
Sobolev inequality holds true:
\begin{equation}\label{loglu3}
\int_{\Omega} v^2 \log \left(\frac{v}{||v||_2}\right)\  |x|^{2-N}
\  dx \le \epsilon \int_{\Omega} |\nabla v|^2 |x|^{2-N}dx +\left
(K_2-\frac{N}{4}\log \epsilon\right) ||v||_2^2\ ,\end{equation}
for any $v\in C^\infty_0(\Omega)$, $v\ge 0$,  and any $\epsilon
>0$; here $K_2$ is a positive constant independent of $\epsilon$ and
$||v||_2:=\left( \int_{\Omega} v^2 |x|^{ 2-N}dx
\right)^\frac{1}{2}$. To see this let us first suppose that the
nonnegative function $v \in C^\infty_0(\Omega)$ is such that
$||v||_2=1$. We then have
$$\int_{\Omega} v^2 \log (v)\  |x|^{2-N} \  dx =\frac{N-2}{4}\int_{\Omega} v^2 \log
\left(v^{\frac{4}{N-2}}\right)\ |x|^{2-N} \  dx \le \frac{N-2}{4}
\log \left(\int_{\Omega} v^{\frac{4}{N-2}+2} \ |x|^{2-N} \ dx
\right)=$$ $$= \frac{N}{4}\log \left(\int_{\Omega}
v^{\frac{2N}{N-2}} \ |x|^{2-N} \ dx \right)^{\frac{N-2}{N}} \le
\frac{N}{4} \log \left (C^{-1}\int_{\Omega}|\nabla v|^2 |x|^{2-N}
dx\right) \ ;$$
 here we have used first Jensen's inequality and then the improved Hardy-Sobolev inequality (\ref{ft}).
For a general nonnegative $v\in C^\infty_0(\Omega)$ we apply the
above inequality to the function $\frac{v}{||v||_2}$, to get
$$\int_{\Omega} v^2 \log \left(\frac{v}{||v||_2}\right) \ |x|^{2-N} \ dx \le \frac{N}{4}
||v||^2_2 \log \left (\frac{C^{-1}}{||v||^2_2}
\int_{\Omega}|\nabla v|^2 |x|^{2-N} dx\right) \ .$$ Since $\log
z\le z$ for any $z>0$, then also $\log y\le \epsilon \frac{4C}{N}y
-\log \left(\epsilon \frac{4C}{N}\right)$, for any $\epsilon>0$;
whence from this we deduce (\ref{loglu3}) with
$K_2:=\frac{N}{4}\log \left(\frac{N}{4C}\right)$.

Inequality (\ref{loglu3}) implies (\ref{logsob2}) via the
following change of variables $u:=v|x|^{\frac{2-N}{2}}$. Finally
from (\ref{logsob1}) and (\ref{logsob2}), the logarithmic Hardy
Sobolev inequality (\ref{logsob3}) easily follows with constant
$K_3:=K_1+K_2+\frac{N+2}{4} \log 2$.

\finedim

We are now ready to give the proof of Theorem \ref{th5}.

\bigskip
{\em Proof of Theorem \ref{th5}: } Let us define, as in Section 2
of \cite{D2}, the operator $\tilde K:=U^{-1} (K -\lambda_1) U$,
$U:L^2(\Omega, \varphi_1^2 dx )\to L^2(\Omega)$ being the unitary
operator $Uw:=\varphi_1 w$, thus $\tilde K:=-\frac{1}{\varphi_1^2}
div (\varphi_1^2 \nabla)$. Here $\varphi_1 >0$ denotes the first
eigenfunction and $\lambda_1>0$ the first eigenvalue corresponding
to the Dirichlet problem $-\Delta \varphi_1
-\frac{(N-2)^2}{4|x|^2} \varphi_1 =\lambda_1 \varphi_1$ in
$\Omega$, $\varphi_1=0$ on $\partial \Omega$, normalized in such a
way that $\int_{\Omega} \varphi_1^2(x)\  dx =1$.  Due to the
results in Lemma 7 in \cite{DD} and using Theorem 7.1 in
\cite{DS1} on one hand and elliptic regularity on the other, there
exist two positive constants $c_1, c_2$ such that

\begin{equation}
\label{auto} c_1 |x|^{\frac{2-N}{2}} d(x)\le \varphi_1(x)\le
c_2|x|^{\frac{2-N}{2}} d(x), ~~~~\forall \ x \in \Omega \ .
\end{equation}
>From this and (\ref{logsob3}) we deduce the following logarithmic
Sobolev inequality \begin{equation}\label{boh} \int_{\Omega} w^2
\log \left(\frac{w}{||w||_2}\right) \  \varphi_1^2 dx \le \epsilon
<\tilde K w, w>_{L^2(\Omega, \varphi_1^2 dx)} +\left(K_4
-\frac{N+2}{4} \log \epsilon \right)||w||^2_2 \ ,
\end{equation}
for any $w\in C^\infty_0(\Omega\setminus \{0\})$, $w\ge 0$, and
any $\epsilon>0$; where $K_4:=K_3 -\log c_1$ and
$||w||_2:=\left(\int_{\Omega} w^2 \varphi_1^2 \
dx\right)^{\frac{1}{2}}$. Let us remark that only the lower bound
in estimate (\ref{auto}) was used.

>From now on one can use the standard approach of \cite{D4} to
complete the proof of the theorem. Here are some details for the
convenience of the reader.

As a first step we claim that the following $L^p$ logarithmic
Sobolev inequalities holds true:

\begin{equation}
\label{loglp} \frac{p}{2}\int_{\Omega} w^p \log
\left(\frac{w}{||w||_p}\right) \ \varphi_1^2 \ dx
 \le \epsilon \ \frac{p}{2} <\tilde K w, w^{p-1}>_{L^2(\Omega, \varphi_1^2 \ dx )} +
 \left(K_4 -\frac{N+2}{4} \log \epsilon \right)||w||^p_p
\end{equation}
for any $w\in C^\infty_0(\Omega\setminus\{0\})$, $w\ge 0$, and any
$\epsilon>0$, $p>2$. To see this we apply inequality (\ref{boh})
to $w^{\frac{p}{2}}$; whence due to the fact that
$$
\int_{\Omega} |\nabla w^\frac{p}{2}|^2  \varphi_1^2\  dx
=\frac{p^2}{4} \int_{\Omega} w^{p-2}|\nabla w|^2 \varphi_1^2\ dx =
\frac{p^2}{4(p-1)} <\nabla w, \nabla w^{p-1}>_{L^2(\Omega,
\varphi_1^2 \ dx)}
 \le \frac{p}{2} <\tilde K w, w^{p-1}>_{L^2(\Omega, \varphi_1^2 \ dx)} \
 ,
$$
since $\frac{p}{2(p-1)}\le 1$ if $p\ge 2$; the claim follows.

Let $H^1_0(\Omega, \varphi_1^2 \ dx)$ be the closure of
$C^\infty_0(\Omega)$ with respect to the norm
$$||w||_{H^1_{0,\varphi_1^2}}:=\left\{\int_{\Omega}
\left(|\nabla w|^2 \ \varphi^2_1 + w^2 \varphi_1^2\right)
dx \right\}^{\frac{1}{2}}\ ; $$ as one can easily prove this is
also the closure of $C^\infty_0(\Omega\setminus \{0\})$ with
respect to the same norm. Then to the operator $\tilde K$ defined
in the domain $D(\tilde K)=\{w \in H^1_0(\Omega, \varphi_1^2 \
dx): \tilde K w \in L^2(\Omega, \varphi_1^2 \ dx )\}$ it is
naturally associated the bilinear symmetric form defined as
follows $\tilde {\mathcal K}[w_1,w_2]:=<\tilde K w_1,
w_2>_{L^2(\Omega, \varphi_1^2 \ dx)}= \int_{\Omega} \nabla w_1
\nabla w_2 \ \varphi_1^2\ dx$,
 which is a Dirichlet form.
Whence Lemma 1.3.4 and Theorems 1.3.2 and 1.3.3 in \cite{D4}
implies that $e^{-\tilde K t}$, which is an analytic contraction
semigroup in $L^2(\Omega, \varphi_1^2 \ dx)$, is also positivity
preserving and a contraction semigroup in $L^p(\Omega, \varphi_1^2
\ dx)$ for any $1\le p\le \infty$. As a consequence for any $t>0$
and any $p\ge 2$  $$e^{-{\tilde K t}}[L^2(\Omega, \varphi_1^2 \
dx)\cap L^\infty (\Omega)]^{+}\subset [H^1_0(\Omega, \varphi_1^2 \
dx)\cap L^p(\Omega, \varphi_1^2 \ dx)\cap L^\infty(\Omega)]^{+}\
;$$ where we denote by $[E]^{+}$ the subset of positive functions
in the space $E$.

\noindent Thus by density argument the $L^p$ logarithmic Sobolev
inequality (\ref{loglp}), more generally applies to any function
in $\cup_{t>0}\ e^{-{\tilde K t}} [L^2(\Omega, \varphi_1^2 dx)\cap
L^\infty (\Omega)]^{+}\ .$

\noindent This means that Theorem 2.2.7 in \cite{D4} can be
applied, in the same way as in Corollary 2.2.8 in \cite{D4}, to
the operator $\tilde K$; whence obtaining that $$||e^{-\tilde K
t}||_{2\to \infty}\le C t^{-\frac{N+2}{4}} \ ,$$ and by duality
that $$||e^{-\tilde K t}||_{1\to 2}\le C t^{-\frac{N+2}{4}} \ ,$$
that is
$$||e^{-\tilde K t}||_{1\to \infty}\le C t^{-\frac{N+2}{2}}\ . $$
Here we use the following notation:
$$||e^{- \tilde K t}||_{q\to p}:=\sup_{0<||f||_q\le 1} \frac{||e^{-\tilde K t} f(x)||_p}{||f(x)||_q}\ ,$$ where
$||f||_q:=\left(\int_{\Omega} |f|^q \varphi_1^2 \ dx
\right)^{\frac{1}{q}} \ .$ This implies, by Dunford-Pettis
theorem, that the semigroup $e^{-\tilde K t}$ is indeed a
semigroup of integral operators; that is a heat kernel $\tilde
k(t,x,y)$ associated to the semigroup  $e^{-\tilde K t}$ is well
defined and satisfies the following pointwise upper bound $\tilde
k(t,x,y)\le C \frac{1}{t} t^{-\frac{N}{2}}$, for any $x,y \in
\Omega$ and any $t>0$. Theorem \ref{th5} then follows, due to the
upper bound in (\ref{auto}) and to the fact that, as a consequence
of the unitary operator $U$, the heat kernels $k(t,x,y)$ and
$\tilde k(t,x,y)$, corresponding respectively to $K$ and $\tilde
K$, satisfy the following equivalence
 \beq \label{stimafi2plus} k(t,x,y)\equiv \varphi_1 (x) \varphi_1 (y) \ \tilde k (t,x,y)
e^{-\lambda_1 t} \ .\eeq

\finedim

\begin{remark}
\label{re1} Applying Davies's method of exponential perturbation
to the operator $\tilde K$ (see Section 2 in \cite{D3} for
details), the upper bound in Theorem \ref{th5} can be improved by
adding a factor $c_\delta e^{-\frac{|x-y|^2}{4(1+\delta)t}}$.
\end{remark}

\medskip

Let us now deduce from the upper bound in Theorem \ref{th5} an
analogous lower bound for time large enough, thus completing the
proof of Theorem \ref{thmoriginbis}. We argue as in Theorem 6 of
\cite{D1} (see also Proposition 4 of \cite{D2}), we give the
details here for the convenience of the reader.

\medskip
{\em Proof of Theorem \ref{thmoriginbis}: } Making use of the same
notation as in the proof of Theorem \ref{th5}, the lower bound we
want to prove corresponds to the statement $\tilde k(t,x,y)\ge C$
for any $x,y\in \Omega$ if $t$ is large enough, $C$ being some
positive constant.

For any $f\in L^1(\Omega, \varphi_1^2 \ dx )$, we clearly have
$$f=<f,1> 1 +g,$$
where $<f,1>:=<f,1>_{L^2(\Omega, \varphi_1^2 \ dx)}$, and
$<g,1>=0$, since $\int_{\Omega} \varphi_1^2(x) dx =1$. Thus,
making
 use of the fact that by definition $\tilde K 1=0$ we have
$$e^{-\tilde Kt} f=<f,1>1+e^{-\tilde K t}g,$$
that is the semigroup $e^{- A t}f :=e^{-\tilde Kt} f-<f,1>1$, to
whom it is clearly associated the heat kernel $\tilde k(t,x,y)-1$,
is such that for any $f\in L^1(\Omega, \varphi_1^2 \ dx)$
$$e^{-A t}f\equiv e^{-\tilde K t}g,$$
where $g=g(f)$ is a function in $L^1(\Omega, \varphi_1^2 dx)$ such
that $<g,1>=0$. Thus, due to Theorem \ref{th5}
$$||e^{- A t}||_{1\to \infty}\le ||e^{-\tilde K t}||_{1\to \infty}\le C t^{-\frac{N+2}{2}},$$
here $C$ is some positive constant; this is equivalent to say that
$$|\tilde k(t,x,y)-1|\le C t^{-\frac{N+2}{2}}\ ,$$
from which the claim easily follows for $t$ large enough.

\finedim

In the sequel we will give the proof of Theorem \ref{prop7}. We
will use the following lemma whose proof will be postponed until
the end of this subsection.

\begin{lemma}
\label{lemma1} For $N\ge 3$, let $\Omega\subset \R^N$ be a smooth
bounded domain containing the origin. Then there exists $\delta_0
>0$, such that
$$\inf_{f\in C^\infty_0(\Omega_\delta)} \frac{\int_{\Omega_\delta} |x|^{2-N} |\nabla f|^2 dx }{\int_{\Omega_\delta}
|x|^{2-N} \frac{f^2}{d^2(x)}dx }= \frac{1}{4} \ ,$$  for all
$0<\delta\le \delta_0$; here $\Omega_\delta:=\{x\in \Omega:
\dist(x,\partial \Omega)\le \delta\}$.
\end{lemma}

We are now ready to prove the improved Hardy inequality.

\medskip
{\em Proof of Theorem \ref{prop7}: } (i) Let us first prove the
claim on any domain $\Omega$ satisfying condition (\ref{cond}). To
this end let us define for any $u\in C^\infty_0(\Omega)$ as a new
variable $w:=|x|^{\frac{N-2}{2}} d^{-\frac{1}{2}}(x) u$, obviously
$w\in H^{1}_0(\Omega)$. By direct computations we have
$$\nabla u= \frac{2-N}{2} |x|^{-\frac{N}{2}-1}\  x \ d^{\frac{1}{2}} \ w + |x|^{\frac{2-N}{2}}\frac {1}{2}
d^{-\frac{1}{2}} \nabla d \ w\  +|x|^{\frac{2-N}{2}} d^{\frac{1}{2}} \nabla w$$
thus
$$\int_{\Omega} |\nabla u|^2 = \int_{\Omega} \left(\frac{(N-2)^2}{4} |x|^{-N} \ d\  w^2 + |x|^{2-N}\frac{1}{4} d^{-1}\
 w^2\  +|x|^{2-N} d \ |\nabla w|^2 \right) dx +$$
$$+\int_{\Omega} \left(-\frac{N-2}{2} |x|^{-N}\  w^2\  x\  \nabla d -(N-2) |x|^{-N} \ d\  w\  x\  \nabla w +\nabla d
\nabla w \ w\  |x|^{2-N}\right) dx \ .$$
Whence
$$\int_{\Omega} \left(|\nabla u|^2-\frac{(N-2)^2}{4|x|^2} u^2 -\frac{1}{4d^2} u^2 \right)dx =
$$
$$=\int_{\Omega}\left(|\nabla w|^2 \ d \  |x|^{2-N} -\frac{N-2}{2} |x|^{-N}\  w^2 \ x \  \nabla d -\frac{(N-2)}{2}
|x|^{-N}\  d \  x\  \nabla w^2 +\frac{1}{2}\  \nabla d \ \nabla w^2\   |x|^{2-N}\right) dx=$$
$$=  \int_{\Omega}\left(|\nabla w|^2 d |x|^{2-N} -\frac{N-2}{2} |x|^{-N} \ w^2 \ x \ \nabla d +\frac{(N-2)}{2}
div (|x|^{-N} d \  x ) w^2 -\frac{1}{2} div (|x|^{2-N} \nabla d) w^2 \right) dx=$$
 $$= \int_{\Omega}\left(|\nabla w|^2 \ d\  |x|^{2-N} -\frac{1}{2} div (|x|^{2-N}\nabla d ) w^2
  \right) dx\ge 0 ,$$
due to condition (\ref{cond}) on $\Omega$. Thus inequality
(\ref{improvedHardy}) is proved with constant $C(\Omega)\equiv
\frac{1}{4}$ in any domain $\Omega$ satisfying condition
(\ref{cond}).

\medskip

(ii) Let us prove indirectly the claim in the remaining case. To
this end let us denote by $H^1_0(\Omega, |x|^{2-N} dx)$ the
closure of $C^\infty_0(\Omega)$ in the norm

\begin{equation}
\label{norma2}
||f||_{H^1_{2-N}}:=\left\{\int_{\Omega} (|\nabla f|^2 +f^2) |x|^{2-N}dx\right \}^{\frac{1}{2}}\ .
\end{equation}
The improved Hardy inequality (\ref{improvedHardy}) we are going
to prove, in the new variable $v:=|x|^{\frac{N-2}{2}} u$ reads as
follows
$$\int_{\Omega}|\nabla v|^2 |x|^{2-N} dx \ge C \int_{\Omega} |x|^{2-N} \frac{v^2}{d^2} dx \ .$$
Let us suppose that the improved Hardy inequality
(\ref{improvedHardy}) is false; whence let us suppose that the
following holds true

$$\inf_{\{\int_\Omega |x|^{2-N} \frac{v^2}{d^2} \ dx \ =\ 1\}}
\int_{\Omega} |x|^{2-N} |\nabla v|^2 dx =0 \ ;$$
\noindent thus there exists a sequence $\{v_j\}_{j\ge 0}$ in
$H^1_0(\Omega, |x|^{2-N} dx)$ such that $\int_\Omega |x|^{2-N}
\frac{v_j^2}{d^2}\  dx =1$, and

\begin{equation}
\label{meno1} \int_{\Omega} |x|^{2-N} |\nabla v_j|^2 dx \to 0 ,
~~~ \hbox { as } j\to \infty \ .
\end{equation}
For any arbitrary function $\varphi \in C^\infty_0(\Omega)$, such
that $\varphi\equiv 1$ in a neighborhood of the origin, we also
have
 $$\int_{\Omega} |x|^{2-N} |\nabla (\varphi v_j)|^2 dx \le 2 \int_{\Omega} |x|^{2-N}
 \left(|\nabla v_j|^2
 \varphi^2 +|\nabla \varphi|^2 v_j^2 \right)dx  \le$$

\begin{equation}
\label{meno2} \le C  \int_{\Omega} |x|^{2-N} \left(|\nabla v_j|^2 + v_j^2 \right)dx \le
C  \int_{\Omega} |x|^{2-N} |\nabla v_j|^2dx  \to 0 \hbox {  as } j\to \infty\ .
\end{equation}
Here we use the fact that the following inequality holds true

\begin{equation}
\label{util} \int_{\Omega} |x|^{2-N} f^2 dx \le C  \int_{\Omega}
|x|^{2-N} |\nabla f|^2 dx,  ~~~ \forall \ f \in H^1_{0}(\Omega,
|x|^{2-N} dx)\ .
\end{equation}
Inequality (\ref{util}) for example follows easily from inequality
(\ref{ft}) by Holder inequality. From estimate (\ref{meno2}) and
inequality (\ref{util}) (applied to $f:=\varphi v_j$)  we easily
deduce that
$$\int_{\Omega} |x|^{2-N} \varphi^2 v_j^2 \to 0, ~~~ \hbox {  as } j\to \infty,$$
or similarly (due to the fact that $\varphi$ has compact support
inside $\Omega$) that

\begin{equation}
\label{zero} \int_{\Omega} |x|^{2-N} \varphi^2 \frac{v_j^2}{d^2}dx
\to 0, ~~~ \hbox {  as } j\to \infty\ .
\end{equation}
We then compute $$1=\int_{\Omega} |x|^{2-N} \frac{v_j^2}{d^2} dx =
\int_{\Omega} |x|^{2-N} \frac{(\varphi v_j+(1-\varphi) v_j)^2}
{d^2}dx = $$
$$=\int_{\Omega} |x|^{2-N} \varphi^2 \frac{v_j^2}{d^2}dx + 2 \int_{\Omega} |x|^{2-N} \varphi (1-\varphi) \frac{v_j^2}
{d^2} dx + \int_{\Omega} |x|^{2-N} (1-\varphi)^2 \frac{v_j^2}{d^2}
dx . $$ We observe that the first two terms in the last line tend
to zero as $j$ tends to infinity and therefore we obtain that

\begin{equation}
\label{stima1}\int_{\Omega} |x|^{2-N} (1-\varphi)^2
\frac{v_j^2}{d^2} dx =1 +o(1), ~~~ \hbox { as } j\to \infty \ .
\end{equation}
On the other hand we have that
$$\int_{\Omega} |x|^{2-N} |\nabla [(1-\varphi)v_j]|^2 dx \le 2 \int_{\Omega} |x|^{2-N} |\nabla v_j|^2 dx  + 2
 \int_{\Omega} |x|^{2-N} |\nabla (\varphi v_j)|^2 dx \ ,$$
both terms in the right hand side going to zero as $j$ tends to infinity due to (\ref{meno1}) and (\ref{meno2});
whence we deduce that

\begin{equation}
\label{stima2} \int_{\Omega} |x|^{2-N} |\nabla [(1-\varphi)v_j]|^2
dx \to 0, ~~~ \hbox { as } j\to \infty \ .
\end{equation}
Since for any $j\ge 0$ the function $f:=(1-\varphi)v_j$ is an
element of $H^1_0(\Omega_\delta)$ for a suitable choice of the
function $\varphi$ (take it identically  one in a subset
containing $\Omega\setminus \Omega_\delta$), by means of
(\ref{stima1}) and (\ref{stima2}) we reach a contradiction with
Lemma \ref{lemma1}, thus proving the improved Hardy inequality.

\finedim

A similar improved Hardy inequality for a potential behaving like
$((N-2)^2/4) |x|^{-2}$ near the origin and exactly like $(1/4)
d^{-2}(x)$ near the boundary is also shown without any geometric
assumption on the domain $\Omega$ (see Theorem \ref{bolo} below).

\bigskip
We next prove Lemma \ref{lemma1}. One can consider it as a
consequence of the following more general
 result.

\begin{lemma}
\label{prop8} For $N\ge 3$, let $\Omega\subset \R^N$ be a smooth
bounded domain.  Then there exists a positive constant
$\delta_0=\delta_0(\Omega)$, such that for any $V\in
L^\frac{N}{2}(\Omega_{\delta_0})$ and $0<\delta \le \delta_0$, we
have the following estimate
$$\int_{\Omega_\delta} \left(|\nabla u|^2 -\frac{1}{4d^2}  u^2 \right) dx \ge
c\int_{\Omega_\delta} Vu^2 dx \ ,~ ~ ~ \forall \ u\in
C^\infty_0(\Omega_\delta)\ ;$$ here $c=c(\delta)\to \infty$ as
$\delta \to 0$ and $\Omega_\delta:=\{x\in \Omega: \dist(x,\partial
\Omega)\le \delta\}$.
\end{lemma}

\medskip
{\em Proof of Lemma \ref{lemma1}:  }  Let us choose $V(x):=
\frac{(N-2)^2}{4|x|^2}$ in Lemma \ref{prop8} above and let us
choose $\delta$ small enough such that $c(\delta)\ge 1$ and
$0<\delta\le \delta_0$, thus we have

\begin{equation}
\label{ineq1}
\int_{\Omega_\delta} \left(|\nabla u|^2 -\frac{1}{4d^2} u^2 \right) dx \ge \int_{\Omega_\delta}
 \frac{(N-2)^2}{4|x|^2}
u^2 dx,
\end{equation}
for any $u\in C^\infty_0(\Omega_\delta)$. For any $f\in
C^\infty_0(\Omega_\delta)$, $u:= f |x|^{\frac{2-N}{2}}$ will be in
$C^\infty_0(\Omega_\delta)$, moreover by easy computations we have
$$\int_{\Omega_\delta} \left(|\nabla u|^2 -\frac{1}{4d^2} u^2 \right) dx = \int_{\Omega_\delta}
\left(|\nabla f|^2 +\frac{(N-2)^2}{4|x|^2} f^2  -\frac{1}{4 d^2} f^2 \right) |x|^{2-N}dx \ ,$$ thus (\ref{ineq1})
can be restated as follows
 $$\int_{\Omega_\delta} \left(|\nabla f|^2 -\frac{1}{4 d^2} f^2 \right) |x|^{2-N}dx \ge 0 \ ;$$
this proves the claim.

\finedim

\bigskip
Whence it only remains to prove Lemma \ref{prop8}. Before doing so
let us observe that inequality (\ref{ineq1}), simply says that the
improved Hardy inequality (\ref{improvedHardy}) indeed holds true
with constant $C(\Omega)=\frac{1}{4}$ whenever the support of the
functions considered is contained in a neighborhood of the
boundary.

\medskip
The proof of Lemma \ref{prop8} makes use of the following improved
Hardy-Sobolev inequality near the boundary stated in Theorem 3 of
\cite{FMT1}, we recall it here for the convenience of the reader:

\begin{proposition}
\label{prop9} For $N\ge 3$, let $\Omega\subset \R^N$ be a smooth
bounded domain. Then there exist positive constants
$\delta_0=\delta_0(\Omega)$ and $C=C(N)$, such that
$$\int_{\Omega_\delta} \left(|\nabla u|^2 -\frac{1}{4 d^2} u^2  \right) dx \ge C \left(\int_{\Omega_\delta}
u^{\frac{2N}{N-2}} dx\right)^{\frac{N-2}{N}} \ , ~~~\forall \ u\in
C^\infty_0(\Omega_\delta), $$ and any $0<\delta \le \delta_0$;
here $\Omega_{\delta}:=\{x\in \Omega: \dist(x,\partial \Omega)\le
\delta\}$.
\end{proposition}

Let us focus here on the fact that in Proposition \ref{prop9} no
convexity assumption on the domain $\Omega$ is made; this is due
to the fact that we only consider functions whose supports are
contained in a neighborhood of the boundary.

\medskip
{\em Proof of the Lemma \ref{prop8}: } By Holder inequality we
have $$\int_{\Omega_\delta} Vu^2 dx \le \left(\int_{\Omega_\delta}
V^\frac{N}{2} dx\right)^{\frac{2}{N}} \left(\int_{\Omega_\delta}
u^{\frac{2N}{N-2}} dx\right)^{\frac{N-2}{N}}\le $$
$$\le \left(\int_{\Omega_\delta} V^\frac{N}{2} dx\right)^{\frac{2}{N}} C(N)^{-1} \int_{\Omega_\delta}
\left(|\nabla u|^2 -\frac{1}{4} \frac{u^2}{d^2} \right) dx\ ,$$
the last step being due to Proposition \ref{prop9}.
This proves the claim with constant $$c(\delta):=\left(\int_{\Omega_\delta}
 V^\frac{N}{2} dx\right)^{-\frac{2}{N}} C(N)\  ,$$ which tends to infinity as $\delta$ tends to zero due to the
 integrability assumption on $V$.

\finedim

\medskip
With some minor changes in the proof of Theorem \ref{prop7} one
can indeed prove the following improved Hardy inequality, which
does not a priori requires the bounded domain $\Omega$ to be
smooth.

\begin{theorem}
\label{prop10} For $N\ge 3$, let $\Omega\subset \R^N$ be a bounded
domain containing the origin such that
$$\int_{\Omega} |\nabla u|^2 dx \ge C \int_{\Omega} \frac{u^2}{d^2} dx, ~~~\forall \ u\in
C^\infty_0(\Omega)$$
and some positive constant $C$. Then
there exists a positive constant $\tilde C$ such that
$$\int_{\Omega} \left(|\nabla u|^2 -\frac{(N-2)^2}{4|x|^2}  u^2\right) dx \ge \tilde C
\int_{\Omega} \frac{u^2}{d^2} dx, ~~~ \forall \ u\in
C^\infty_0(\Omega)\ .$$

\end{theorem}

We finally mention the following related new Hardy inequality,
which we think is of independent interest

\begin{theorem}
\label{bolo} For $N\ge 3$, let $\Omega\subset \R^N$ be a smooth
bounded domain containing the origin, and define for $\epsilon>0$,
$$V_\epsilon(x)=\begin{cases} \frac{(N-2)^2}{4|x|^2} & if \ \{x \in \Omega: d(x)\ge \epsilon \} \\
\frac{1}{4d^2(x)} & if \ \{x\in \Omega: d(x)<\epsilon\} \
.\end{cases}
$$ Then there exists $\epsilon_0=\epsilon_0(\Omega)$ such that for
all $0<\epsilon\le \epsilon_0$ and $u\in C^\infty_0(\Omega)$, we
have
$$\int_{\Omega} |\nabla u|^2 dx \ge \int_{\Omega} V_\epsilon(x) u^2  dx \ .$$
\end{theorem}

\smallskip

{\em Proof: } We will only sketch it. Let $\Omega_1=\{x\in \Omega:
d(x)\ge  \epsilon\}$. Then using  the change of variable
$u:=|x|^\frac{2-N}{2} v$, one can prove the following inequality
$$\int_{\Omega_1} \left(|\nabla u|^2 -
\frac{(N-2)^2}{4|x|^2} u^2 \right) dx \ge
\frac{2-N}{2}\int_{\partial \Omega_1} \frac{u^2}{|x|^2} x \cdot
\nu \ dS_x. $$ Similarly using the change of variable
$u:=d^\frac{1}{2}(x) X^{-\frac{1}{2}}(d(x)) v$ with $X(t)=(1-\ln
t)^{-1}$ one can prove the following inequality
$$\int_{\Omega\setminus \Omega_1} \left(|\nabla u|^2 -
\frac{1}{4d^2(x)} u^2 \right) dx \ge - \ \frac{1}{4}\int_{\partial
\Omega_1} \frac{u^2}{d(x)} \nabla d \cdot \nu \ dS_x, $$  for any
$0<\epsilon\le \min\{e^{-1}, \epsilon_1\}$ where $\epsilon_1>0$ is
such that $d^{-1}X(d)+2\Delta d(\ln d)\ge 0$ for $d\le
\epsilon_1$. The result then follows showing that for $0<\epsilon
\le \epsilon_0= \min \{e^{-1}, \epsilon_1, \frac{R}{2N-3}\}$ we
have $\left[\frac{2(2-N)}{|x|^2} x -\frac{1}{d(x)} \nabla d\right]
\cdot \nu \ge 0$ since $\nu:=-\nabla d$ on $\partial \Omega_1$;
here $R$ denotes a positive constant such that $B(0,R)\subset
\Omega$, which exists due to the assumption on $\Omega$. \finedim

\subsection{Complete sharp description of the heat kernel for small values of time}

In this section we prove the two-sided sharp estimate on the heat
kernel $k(t,x,y)$ stated for small time in Theorem
\ref{thmorigin}.

\medskip

{\em Proof of Theorem \ref{thmorigin}} Since for any $x\in \Omega$
and for some positive constants $c_1,c_2$ we have the following
estimate $c_1 |x|^{\frac{\lambda}{2}} d^{\frac{\alpha}{2}}(x) \le
\varphi_1(x)\le c_2 |x|^{\frac{\lambda}{2}}
d^{\frac{\alpha}{2}}(x)$ for $\alpha=2$ and $\lambda=2-N$, we can
apply the result of Theorem \ref{stab} to the operator $\tilde
K=-\frac{1}{\varphi_1^2(x)}div(\varphi_1^2(x) \nabla)$. Hence due
to (\ref{stimafi2plus}) the result follows immediately.

\finedim

Let us finally make some remarks concerning Schr\"odinger
operators having potential $V(x)=c|x|^{-2}$. Arguing as in Lemma 7
in \cite{DD} one can easily prove that the first Dirichlet
eigenfunction for the Schr\"odinger operator $-\Delta
-\frac{c}{|x|^2} $, $0<c<\frac{(N-2)^2}{4}$, behaves like
$|x|^{\frac{\lambda}{2}} d(x)$ on all $\Omega$, where
$\lambda:=2-N+\sqrt{(N-2)^2-4c}$. Then we have

\begin{theorem}
For $N\ge 3$, let $\Omega \subset \R^N$ be a smooth bounded domain
containing the origin. Then there exist positive constants $C_1,
C_2$, with $C_1\le C_2$, and $T>0$ depending on $\Omega$ such that
$$C_1 \min\Big \{(|x|+\sqrt{t})^{\frac{|\lambda|}{2}}(|y|+\sqrt{t})^{\frac{|\lambda|}{2}}  ,
\frac{d(x) d(y)}{t} \Big\}   |xy|^{\frac{\lambda}{2}} t^{-\frac{N}{2}} e^{-C_2 \frac{|x-y|^2}{t}} \le $$
$$\le k_c(t,x,y)\le C_2 \min\Big \{(|x|+\sqrt{t})^{\frac{|\lambda|}{2}}(|y|+\sqrt{t})^{\frac{|\lambda|}{2}},
\frac{d (x) d(y)}{t} \Big\}  |xy|^{\frac{\lambda}{2}}
t^{-\frac{N}{2}} e^{-C_1 \frac{|x-y|^2}{t}} \ ,$$ for all $x,y\in
\Omega$ and $0<t\le T$; here $k_c(t,x,y)$ denotes the heat kernel
associated to the operator $-\Delta -\frac{c}{|x|^2}$ in $\Omega$
under Dirichlet boundary conditions for $0<c<\frac{(N-2)^2}{4}$,
and $\lambda:=2-N+\sqrt{(N-2)^2-4c}$.
\end{theorem}

\begin{theorem}
For $N\ge 3$, let $\Omega \subset \R^N$ be a smooth bounded domain
containing the origin. Then there exist two positive constants
$C_1, C_2$, with $C_1\le C_2$, such that
$$C_1 \ d(x)\  d(y) \ |xy|^{\frac{\lambda}{2}} e^{-\lambda_1 t}
\le k_c(t,x,y)\le C_2\  d (x) \ d(y) \ |xy|^{\frac{\lambda}{2}}
e^{-\lambda_1 t} \ ,$$ for all $x,y\in \Omega$ and $t>0$ large
enough; here $k_c(t,x,y)$ denotes the heat kernel associated to
the operator $-\Delta -\frac{c}{|x|^2}$ in $\Omega$ under
Dirichlet boundary conditions for  $0<c<\frac{(N-2)^2}{4}$,
$\lambda_1$ its (positive) elliptic first eigenvalue and
$\lambda:=2-N+\sqrt{(N-2)^2-4c}$.
\end{theorem}

\smallskip

\medskip

\setcounter{equation}{0}
\section{Critical boundary singularity}

In this section we prove Theorems \ref{crit} and \ref{thmbrybis}
as well as a new Hardy-Moser inequality (Theorem \ref{prop4}). The
structure of this section is as follows.

In Subsection 4.1 we first prove the improved Hardy-Moser
inequality.  Then in Subsection 4.2 we get the global in time
pointwise upper bound for the heal kernel of the Schr\"odinger
operator $-\Delta -(1/4)d^{-2}(x)$, which is sharp when $x$ and
$y$ are close to the boundary (see Theorem \ref{th3}). Then
arguing as in \cite{D1}, we deduce the sharp heat kernel lower
bound for time large enough, thus completing the proof of Theorem
\ref{thmbrybis}.

The proof of Theorem \ref{crit} is finally completed in Subsection
4.3, using the parabolic Harnack inequality up to the boundary
stated in Theorem \ref{harnack}.

\subsection{The improved Hardy-Moser inequality}

Here we will prove a new  improved Hardy-Moser inequality
which we think it is of independent interest.
 The proof is based on an auxiliary Hardy-Sobolev
inequality, that we will show here, as well as on the following
improved Hardy inequality stated in Theorem A in \cite{BFT1}.

\begin{proposition}
\label{prop3} For $N\ge 2$, let $\Omega  \subset  \R^N$ be a
smooth bounded and convex domain. Then there exists $D_0$ positive
such that for all $D\ge D_0$
$$\int_{\Omega}\left( |\nabla u|^2 - \frac{1}{4d^2(x)} u^2 \right) \  dx \ge \frac{1}{4} \int_{\Omega}
\frac{X^2\left(\frac{d(x)}{D}\right)}{d^2(x)}  u^2 \  dx \  ,$$
for any $u \in C^\infty_0(\Omega)$;  here $X(t):=\frac{1}{1-\ln
t}$, $t\in (0,1]$.
\end{proposition}

\medskip

Let us now state the auxiliary Hardy-Sobolev inequality we will
use in the sequel.

\begin{lemma}
\label{prop2} Let $\alpha>0$, $N\ge 2$ and $\Omega  \subset \R^N$
be a smooth bounded domain. Then there exist $\delta_0
>0$ and $C=C(\alpha,\delta_0)>0$ such that
$$\int_{\Omega}d^\alpha(x)|\nabla v|  \ dx  + \int_{\Omega\setminus \Omega_\delta} d^{\alpha-1}(x)|v|  \  dx \ge
C \left( \int_{\Omega}  d^{\frac{\alpha N}{N-1}}(x)
|v|^{\frac{N}{N-1}}\  dx \right)^{\frac{N-1}{N}}\ ,$$ for any $v
\in C^\infty_0(\Omega)$ and any $0<\delta \le \delta_0$;  here
$\Omega_\delta:=\{x\in \Omega : \dist(x,\partial \Omega)\le
\delta\}$.
\end{lemma}

\smallskip

{\em Proof: } We will follow closely the argument of \cite{FMT2}.
Our starting point is the following Gagliardo-Nirenberg inequality
(see p. 189 in \cite{M})
$$S_N||f||_{L^{\frac{N}{N-1}}(\Omega)}\le ||\nabla f||_{L^1(\Omega)}\ , ~ ~ ~ \forall \
f\in C^\infty_0(\Omega) \ ,$$ where $S_N$ is a positive constant
depending only on $N$.

For any $v \in C^\infty_0(\Omega)$ let us apply the above
inequality to the function $f:=d^\alpha v$. Hence we obtain
$$S_N ||d^\alpha v||_{L^{\frac{N}{N-1}}(\Omega)}\le \int_{\Omega} d^\alpha(x)|\nabla v| dx +
\alpha \int_{\Omega} d^{\alpha-1}(x) |v| dx \ .$$ We next estimate
the last term above. Let $\varphi_\delta \in
C^\infty_0(\Omega_{2\delta})$, $0\le \varphi_\delta \le 1$, be a
cut off function which is identically one in $\Omega_{\delta}$ and
identically zero in $\R^N\setminus \Omega_{2\delta}$. Clearly
$v=\varphi_\delta v+(1-\varphi_{\delta})v$. Then we have
$$\alpha \int_{\Omega} d^{\alpha-1}(x) |v| dx\le
 \alpha \int_{\Omega} d^{\alpha-1}(x) |\varphi_{\delta}v| dx + \alpha
\int_{\Omega} d^{\alpha-1}(x) (1-\varphi_\delta) |v| dx \le
$$
$$\le\alpha \int_{\Omega} d^{\alpha-1}(x) |\varphi_{\delta}v| dx +
\alpha \int_{\Omega\setminus \Omega_{\delta}} d^{\alpha-1}(x) |v|
dx \ .$$ Concerning the first term on the right hand side we have
$$\alpha \int_{\Omega} d^{\alpha-1}(x) |\varphi_{\delta}v| dx = \int_{\Omega}  \nabla
d^{\alpha}\cdot
 \nabla d |\varphi_{\delta}v| dx = -\int_{\Omega} d^\alpha(x) \nabla d\cdot \nabla |\varphi_{\delta}v|
 dx -
 \int_{\Omega} d^\alpha(x)\Delta d |\varphi_{\delta}v|dx \le$$
 $$\le
 \int_{\Omega} d^\alpha(x) |\nabla (\varphi_{\delta}v)| dx +c_0
 \delta \int_{\Omega}d^{\alpha-1}(x)|\varphi_{\delta}v| dx \ ,$$
 here we used the smoothness assumption on $\Omega$ which
 implies that $|d\Delta d|\le c_0\delta$ in $\Omega_\delta$ for
 $\delta$ small, say $0<\delta\le \delta_0$, and for some positive constant $c_0$ independent
 of $\delta$ ($\delta_0$, $c_0$ depending on $\Omega$). Thus we have for any $0<\delta \le \delta_0$
$$\alpha \int_{\Omega} d^{\alpha-1}(x) |\varphi_{\delta}v| dx \le
\frac{\alpha}{\alpha-c_0\delta_0}\int_{\Omega} d^\alpha(x) |\nabla
(\varphi_{\delta}v)| dx\le C\int_{\Omega} d^\alpha(x) |\nabla v|
\varphi_{\delta} dx+ $$ $$+\frac{C}{\delta}
\int_{\Omega_{2\delta}\setminus \Omega_{\delta}} d^\alpha(x) |v|
dx\le C\int_{\Omega} d^\alpha(x) |\nabla v| dx+
C\int_{\Omega_{2\delta}\setminus \Omega_{\delta}} d^{\alpha-1}(x)
|v|dx \ ,$$ from which the result follows.

\finedim

\medskip
\noindent We next state the  new improved  Hardy-Moser
inequality.

\begin{theorem}
{\bf (Improved Hardy-Moser inequality)} \label{prop4} For $N\ge
2$, let $\Omega \subset  \R^N$ be a smooth bounded and convex
domain. Then there exists a positive constant $C$ such that
$$\left\{\int_{\Omega}\left( |\nabla u|^2 - \frac{1}{4 d^2} u^2 \right) \  dx \right\}
\left(\int_{\Omega} u^2 dx\right )^{\frac{2}{N}}\ge C
\int_{\Omega} u^{2(1+\frac{2}{N})}  dx \ , ~ ~ ~ \forall \ u \in
C^\infty_0(\Omega)\ .$$
\end{theorem}
{\em Proof: } Changing variables by $v:=u
d^{-\frac{1}{2}}$, we get
$$\int_{\Omega} u^{\frac{2(N+2)}{N}} dx  = \int_{\Omega}d^\frac{N+2}{N} v^\frac{2(N+2)}{N}  dx =
\int_{\Omega} d^\frac{\alpha N}{N-1}(v^{2\alpha})^\frac{N}{N-1} dx
\ ,$$ with $\alpha:=\frac{(N+2)(N-1)}{N^2}$. Applying Lemma
\ref{prop2} to the function $v^{2\alpha}$ we have
$$\int_{\Omega} d^\frac{\alpha N}{N-1}(v^{2\alpha})^\frac{N}{N-1}  dx\le
 C \left( \int_\Omega d^\alpha |\nabla v^{2\alpha}| dx +
 \int_{\Omega\setminus\Omega_{\delta}} d^{\alpha-1} v^{2\alpha} dx \right)^{\frac{N}{N-1}}\le $$
$$\le C \left( 2 \alpha \int_\Omega d^\alpha|\nabla v| |v|^{2\alpha-1} dx +
\int_{\Omega\setminus\Omega_{\delta}} d^{\alpha-1} v^{2\alpha} dx
\right)^{\frac{N}{N-1}}\le$$
$$\le C \left\{\left(\int_\Omega d(x)|\nabla v|^2 \ dx\right)^\frac{1}{2}
\left(\int_{\Omega}d^{2\alpha-1}v^{2(2\alpha-1)} dx
\right)^\frac{1}{2}+\left(\int_{\Omega\setminus\Omega_\delta}
\frac{v^2}{d} dx
\right)^\frac{1}{2}\left(\int_{\Omega\setminus\Omega_{\delta}}
d^{2\alpha-1} v^{2(2\alpha-1)} dx \right)^{\frac{1}{2}}
\right\}^{\frac{N}{N-1}}\le$$
$$\le C
\left(\int_{\Omega} d^{2\alpha-1}v^{2(2\alpha-1)} dx
\right)^\frac{N}{2(N-1)} \left\{\int_\Omega \left(|\nabla v|^2 d
-\frac{1}{2}\Delta d \ v^2 \right) dx \right\}^\frac{N}{2(N-1)} \
;
$$
here we used Proposition \ref{prop3} (observe that $\frac{1}{4}
X^2(\frac{d}{D})\ge \frac{1}{4}X^2(\frac{\delta}{D})$ if $x\in
\Omega\setminus \Omega_\delta$) and standard estimates. Returning
to the original variable $u$, we obtain
$$\int_{\Omega} u^{\frac{2(N+2)}{N}} dx \le C
\left(\int_{\Omega} u^{2(2\alpha-1)} dx \right)^\frac{N}{2(N-1)}
\left\{\int_\Omega \left(|\nabla u|^2  -\frac{1}{4d^2} u^2\right)
dx  \right\}^\frac{N}{2(N-1)} \ ,
$$
that is,
$$\left(\int_{\Omega} u^{\frac{2(N+2)}{N}} dx \right)^{\frac{2(N-1)}{N}} \le
 C \left(\int_{\Omega} u^{2(2\alpha-1)} dx \right)\left\{\int_\Omega \left(|\nabla u|^2
  -\frac{1}{4d^2} u^2\right) dx  \right\} \ .
$$
If $N=2$ we have that $\alpha=1$, thus the above inequality
becomes
$$\int_{\Omega} u^{4} dx  \le C \left(\int_{\Omega} u^{2} dx
 \right)\left\{\int_\Omega \left(|\nabla u|^2  -\frac{1}{4d^2} u^2\right) dx  \right\}
$$
which is the sought for estimate. For $N\ge 3$, we use H\"older
inequality  to obtain
$$\left(\int_{\Omega} u^{\frac{2(N+2)}{N}} dx \right)^{\frac{2(N-1)}{N}}
 \le C \left(\int_{\Omega} u^{2} dx \right)^{\frac{2}{N}}
 \left(\int_{\Omega} u^{\frac{2(N+2)}{N}} dx \right)^{\frac{N-2}{N}}
 \left\{\int_\Omega \left(|\nabla u|^2  -\frac{1}{4d^2} u^2\right) dx  \right\}
$$
from which
$$\int_{\Omega} u^{\frac{2(N+2)}{N}} dx
\le C \left(\int_{\Omega} u^{2} dx \right)^{\frac{2}{N}}
 \left\{\int_\Omega \left(|\nabla u|^2  -\frac{1}{4d^2} u^2\right) dx  \right\} \ ;
$$
and this completes the proof of Theorem \ref{prop4}.

\finedim

\subsection{Boundary upper bounds and complete sharp description of the heat kernel
for large values of time}

Here we will first prove the following:
\begin{theorem}
\label{th3} For $N\ge2$, let $\Omega \subset \R^N$ be a smooth
bounded and convex domain. Then there exists a positive constant
$C$ such that
$$h(t,x,y)\le C \frac{d^{\frac{1}{2}}(x) d^{\frac{1}{2}}(y)}{t^\frac{1}{2}}
t^{-\frac{N}{2}}, ~ ~ ~ \forall \ x,y\in \Omega , \ t>0 \ .$$
\end{theorem}
To this end we need the following estimate of \cite{FMT2}:
\begin{proposition}
\label{prop5} For $N\ge2$, let $\Omega \subset \R^N$ be a smooth
bounded and convex domain. Then there exists a positive constant
$C$ such that \begin{equation} \label{sky2} \int_{\Omega}
\left(|\nabla u|^2 -\frac{1}{4 d^2(x)} u^2 \right)dx
 \ge C \left(\int_{\Omega}  d^{\frac{q}{2}(N-2) -N}(x) |u|^{q}\ dx
 \right)^\frac{2}{q} \ ,
\end{equation} for any $u\in C^\infty_0(\Omega)$ and any $2<q\le
\frac{2N}{N-2}$ if $N\ge 3$ or any $2<q<\infty$ if $N=2$.
\end{proposition}

Using (\ref{sky2}) the following logarithmic Sobolev inequality
can be easily obtained
\begin{equation}
\label{log3} \int_{\Omega} v^2 \log \left(\frac{v}{||v||_2}\right)
d \ dx \le \epsilon \int_{\Omega} \left(|\nabla v|^2 d-\frac{1}{2}
\Delta d\  v^2\right) dx +\left(K_1-\frac{N+1}{4} \log
\epsilon\right) ||v||^2_2 \ , \end{equation} for all $v\in
C^\infty_0(\Omega)$, $v\ge 0$, and any $\epsilon >0$; here $K_1$
is a positive constant independent of $\epsilon$ and
$||v||_2:=\left(\int_\Omega |v|^2 \ d \ dx \right)^\frac{1}{2}$.

To obtain (\ref{log3}) we apply (\ref{sky2}) to $v:=u
d^{-\frac{1}{2}}$ to get for any $v\in C^\infty_0(\Omega)$
$$
\int_{\Omega} \left(|\nabla v|^2 d -\frac{1}{2} \Delta d\
v^2\right)dx \ge C \left(\int_{\Omega} v^{q}d^{\frac{q}{2}(N-2)-N
+\frac{q}{2}} \ dx \right)^\frac{2}{q} \ .
$$
Taking $q:=\frac{2(N+1)}{N-1}$ we have
\begin{equation}
\label{sob1}
\int_{\Omega} \left(|\nabla v|^2 d -\frac{1}{2} \Delta d\  v^2\right)dx \ge
C \left(\int_{\Omega}  v^{\frac{2(N+1)}{(N-1)}}d \ dx  \right)^\frac{N-1}{N+1} \ .
\end{equation}
Then arguing in a quite similar way as in the proof of
(\ref{loglu3}) in Subsection 3.1 we obtain (\ref{log3}) with
$K_1:=\frac{N+1}{4}\log\left(\frac{N+1}{4C}\right)$.

{\em Proof of Theorem \ref{th3}: } Let $H^1_0(\Omega, d \ dx)$ be
the closure of $C^\infty_0(\Omega)$ with respect to the norm
$$||v||_{H^1_{0,d}}:=\left\{\int_{\Omega} \left(|\nabla v|^2 \ d +\frac{1}{2}
 (-\Delta d)v^2 \right) dx \right\}^{\frac{1}{2}}\ . $$ Let
$\bar H:=U^{-1} H U$, $U:L^2(\Omega, d \ dx)\to L^2(\Omega)$ being
the unitary operator $Uv:=d^{\frac{1}{2}} v$, thus
 $\bar H:=-\frac{1}{d} div (d \nabla )-\frac{1}{2}\frac{\Delta d}{d}$.
To the operator $\bar H$ defined in the domain $D(\bar H)=\{v \in
H^1_0(\Omega, d \ dx): \bar H v \in L^2(\Omega, d \ dx )\}$ it is
naturally associated the bilinear symmetric form defined as
follows $\bar{\mathcal H}[v_1,v_2]:=<\bar H v_1, v_2>_{L^2(\Omega,
d \ dx)}= <v_1,v_2>_{H^1_0(\Omega, d \ dx)}$  which is a Dirichlet
form. Whence Lemma 1.3.4 and Theorems 1.3.2 and 1.3.3 in \cite{D4}
implies that $e^{-\bar H t}$, which is an analytic contraction
semigroup in $L^2(\Omega, d \ dx)$, is also positivity preserving
and a contraction semigroup in $L^p(\Omega, d \ dx)$ for any $1\le
p\le \infty$. As a consequence for any $t>0$ and any $p\ge 2$
$$e^{-{\bar H t}}[L^2(\Omega, d \ dx)\cap L^\infty
(\Omega)]^{+}\subset [H^1_0(\Omega, d \ dx)\cap L^p(\Omega, d \
dx)\cap L^\infty(\Omega)]^{+}\ ;$$ thus by density argument the
$L^p$ logarithmic Sobolev inequality, which can be deduced as
usual from the $L^2$ logarithmic Sobolev inequality (\ref{log3})
(see Subsection 3.1 where a similar argument is used) more
generally applies to any function in $\cup_{t>0}\ e^{-{\bar H
t}}[L^2(\Omega, d \ dx)\cap L^\infty (\Omega)]^{+}$. This means
that Theorem 2.2.7 in \cite{D4} can be applied, as in Corollary
2.2.8 in \cite{D4},
 to the operator $\bar H$; whence obtaining that
 $$||e^{-\bar H t}||_{2\to \infty}\le C t^{-\frac{N+1}{4}}\ ,$$ and by duality that
 $$||e^{-\bar H t}||_{1\to 2}\le C t^{-\frac{N+1}{4}} \ ,$$ that
 is
 $$||e^{-\bar H t}||_{1\to \infty}\le C t^{-\frac{N+1}{2}}\ .$$
Here we use the following notation:
$$||e^{-\bar Ht}||_{q\to p}:=\sup_{0<||f||_q\le 1} \frac{||e^{-\bar Ht} f(x)||_p}{||f(x)||_q}\ ,$$
where $||f||_q:=\left(\int_{\Omega} |f|^q \ d \ dx
\right)^{\frac{1}{q}}$. This implies, by Dunford-Pettis theorem,
that the semigroup $e^{-\bar H t}$ is indeed a semigroup of
integral operators; that is a heat kernel $\bar h(t,x,y)$
associated to the semigroup $e^{-\bar H t}$ is well defined and
satisfies the following pointwise upper bound $\bar h(t,x,y)\le C
\frac{1}{t^{\frac{1}{2}}} t^{-\frac{N}{2}}$, for any $x,y\in
\Omega$ and any $t>0$. Theorem \ref{th3} then follows, due to the
fact that, as a consequence of the unitary operator $U$, the heat
kernels $h(t,x,y)$ and $\bar h(t,x,y)$, corresponding respectively
to $H$ and $\bar H$, satisfy the following equivalence
$h(t,x,y)\equiv d^{\frac{1}{2}} (x) d^{\frac{1}{2}}(y) \ \bar h
(t,x,y) \ .$

\medskip

\begin{remark}
Applying Davies's method of exponential perturbation to the
operator $\bar H$  (see Section 2 in \cite{D3} for details) the
upper bound in Theorem \ref{th3} can be improved by adding a
factor $c_\delta e^{-\frac{|x-y|^2}{4(1+\delta)t}}$.
\end{remark}

\finedim

\medskip

Let us now give the sketch of the proof of Theorem \ref{thmbrybis}.

\medskip
{\em Proof of Theorem \ref{thmbrybis}: } Let us first improve by
an exponential decreasing in time factor the global in time upper
bound stated in Theorem \ref{th3}. To this end let us define, as
in Section 2 of \cite{D2}, the operator $\tilde
H:=U^{-1}(H-\lambda_1)U$,
 $U:L^2(\Omega, \varphi_1^2 dx )\to L^2(\Omega)$ being  the unitary operator
  $Uw:=\varphi_1 w$, thus $\tilde H:=-\frac{1}{\varphi_1^2} div (\varphi_1^2 \nabla)$.
  Here $\varphi_1>0$ denotes the first eigenfunction and
$\lambda_1>0$ the first
 eigenvalue corresponding to the Dirichlet problem $-\Delta \varphi_1-\frac{1}{4d^2(x)}
  \varphi_1=\lambda_1 \varphi_1$ in $\Omega$, $\varphi_1=0$ on $\partial \Omega$,
  normalized in such a way that $\int_{\Omega} \varphi_1^2(x)\  dx =1$.  Since it is known
  that there exist two positive constants $c_1,c_2$ such that
\begin{equation}\label{stimafi}
  c_1 \ d^{\frac{1}{2}}(x)\le
  \varphi_1(x) \le c_2 \ d^{\frac{1}{2}}(x), ~ \forall \ x \in
  \Omega\end{equation}
  (as a consequence of Lemma 7 in \cite{DD}), logarithmic Sobolev inequalities analogous
   to (\ref{log3}) also hold true if we replace $\bar H$ by $\tilde H$. Thus as a consequence
    of Gross theorem, exactly as in the proof of Theorem \ref{th3}, the corresponding heat
     kernel  $\tilde h(t,x,y)$ satisfies the same pointwise upper bound as $\bar h(t,x,y)$ that
      is $\tilde h(t,x,y)\le \frac{C}{t^\frac{1}{2}}t^{-\frac{N}{2}}$ for any $x,y\in \Omega$
      and any $t>0$. From the definition of $U$, it follows \begin{equation}
      \label{stimafi2} h(t,x,y)\equiv \varphi_1(x)
      \varphi_1(y) \tilde h(t,x,y)e^{-\lambda_1 t}\ ,\end{equation} thus we get that $h(t,x,y)\le C
      \frac{d^\frac{1}{2}(x)d^\frac{1}{2}(y)}{t^{\frac{1}{2}}} t^{-\frac{N}{2}}
      e^{-\lambda_1 t}$ for any $x,y\in \Omega$ and any $t>0$.
       Finally arguing as in Theorem 6 of \cite{D1},
       an analogous lower estimate can be easily deduced (we refer to the proof of
        Theorem \ref{thmoriginbis} where a similar argument is used).

\finedim

\subsection{Complete sharp description of the heat kernel for small values of time}

In this section we prove the two-sided sharp estimate on the heat
kernel $h(t,x,y)$ stated for small time in Theorem \ref{crit}.

Before doing so let observe that Theorem \ref{harnack} entails
also the following parabolic Harnack inequality for the
Schr\"odinger operator having critical singularity at the boundary

\begin{theorem}
\label{phu} For $N\ge 2$, let $\Omega\subset \R^N$ be a smooth
bounded and convex domain. Then there exist positive constants
$C_H$ and $R=R(\Omega)$ such that for $x\in \Omega$, $0<r<R$ and
for any positive solution $u(y,t)$ of $\frac{\partial u}{\partial
t}=\Delta u+\frac{1}{4d^2(y)}u$ in $\left\{\mathcal B(x,r)\cap
\Omega\right\}\times (0,r^2)$ we have the estimate
$${\it ess~sup}_{(y,t)\in \left\{\mathcal B(x,\frac{r}{2})\cap \Omega \right\}\times
(\frac{r^2}{4},\frac{r^2}{2})} \ \  u(y,t) d^{-\frac{1}{2}}(y)\le
C_H ~ {\it ess~inf}_{(y,t)\in \left\{\mathcal B(x,\frac{r}{2})\cap
\Omega \right\}\times (\frac{3}{4} r^2,r^2)} \ \ u(y,t)
d^{-\frac{1}{2}}(y)\ .$$
\end{theorem}

{\em Proof: } As a first step let us observe that if $u$ satisfies
$u_t=-Hu$ then $v(y,t):=e^{\lambda_1 t} \varphi_1(y)^{-1} u(y,t)$
satisfies $v_t=-\tilde H v$. Whence  as a consequence of
(\ref{stimafi}), due to Remark \ref{stab2}, $v$ satisfies Theorem
\ref{harnack} for $\alpha=1$. From this the claim can be easily
deduced.

\finedim

The proof of Corollary \ref{bbb} is similar to the above proof of
Theorem \ref{phu}, thus we omit the details.

\medskip

We are now ready to prove Theorem \ref{crit}.

\smallskip

{\em Proof of Theorem \ref{crit}} Since for any $x\in \Omega$ and
for some positive constants $c_1,c_2$ we have the following
estimate $c_1 d^{\frac{\alpha}{2}}(x) \le \varphi_1(x)\le c_2
d^{\frac{\alpha}{2}}(x)$ for $\alpha=1$, we can apply the result
of Theorem \ref{stab} to the operator $\tilde
H=-\frac{1}{\varphi_1^2(x)}div(\varphi_1^2(x) \nabla)$. Hence due
to (\ref{stimafi2}) the result follows immediately.

\finedim

The proof of Corollary \ref{coint} is similar to the above proofs
of Theorems \ref{crit} and \ref{thmbrybis}, thus we omit the
details.

\bigskip

Let us finally make some remarks concerning Schr\"odinger
operators having potential $V(x)=cd^{-2}(x)$. Arguing as in Lemma
7 in \cite{DD} one can easily prove that the first Dirichlet
eigenfunction for the Schr\"odinger operator $-\Delta
-\frac{c}{d^2(x)} $, $0<c<\frac{1}{4}$, behaves like
$d^\frac{\alpha}{2}$ on all $\Omega$, for $\alpha
:=1+\sqrt{1-4c}$. Then we have:

\begin{theorem}
\label{sopra} For $N\ge 2$, let $\Omega\subset \R^N$ be a smooth
bounded and convex domain. Then there exist positive constants
$C_1,C_2$, with $C_1\le C_2$, and $T>0$ depending on $\Omega$ such
that
$$C_1 \min\left\{1,\frac{d^{\frac{\alpha}{2}}(x)d^{\frac{\alpha}{2}}(y)}{t^{\frac{\alpha}{2}}}\right\}
t^{-\frac{N}{2}} e^{-C_2\frac{|x-y|^2}{t}}\le h_c(t,x,y)
 \le C_2 \min\left\{1,\frac{d^{\frac{\alpha}{2}}(x) d^{\frac{\alpha}{2}}(y)}{t^{\frac{\alpha}{2}}}
 \right\} t^{-\frac{N}{2}} e^{-C_1\frac{|x-y|^2}{t}}$$ for all
$x,y\in \Omega$ and $0<t\le T$; where $h_c(t,x,y)$ denotes the
heat kernel associated to the operator $-\Delta -\frac{c}{d^2(x)}$
in $\Omega$ under Dirichlet boundary conditions, for any
$0<c<\frac{1}{4}$ and $\alpha:=1+\sqrt{1-4c}$.
\end{theorem}

\begin{theorem}
For $N\ge 2$, let $\Omega \subset \R^N$ be a smooth bounded and
convex domain. Then there exist two positive constants $C_1, C_2$,
with $C_1\le C_2$, such that
$$C_1 \ d^{\frac{\alpha}{2}} (x) \ d^{\frac{\alpha}{2}}(y) \ e^{-\lambda_1 t} \le  h_c(t,x,y)\le
 C_2 \ d^{\frac{\alpha}{2}} (x) \ d^{\frac{\alpha}{2}}(y) \ e^{- \lambda_1 t} $$
for all $x,y \in \Omega$ and $t>0$ large enough; where
$h_c(t,x,y)$ denotes the heat kernel associated to the operator
$-\Delta -\frac{c}{d^2(x)}$ in $\Omega$ under Dirichlet boundary
conditions, $\lambda_1$ its (positive) elliptic first eigenvalue,
for any $0<c<\frac{1}{4}$ and $\alpha:=1+\sqrt{1-4c}$.
\end{theorem}

\begin{remark}
\label{noconv} Let us at this point remark that Theorems
\ref{crit} and \ref{sopra} as well as Theorem \ref{phu} concerning
respectively sharp asymptotic for small time and the parabolic
Harnack inequality up to the boundary for the Schr\"odinger
operator having potential $V(x)=cd^{-2}(x)$, hold true also
without any convexity assumption on the domain $\Omega$ under
consideration.
\end{remark}

\setcounter{equation}{0}
\subsection{On Davies conjecture}

In this subsection we consider Davies conjecture. For this we
suppose that $\widetilde{E}$ denotes the self-adjoint operator
associated with the closure of the positive quadratic form
$$Q(f)=\int_{\Omega} \left(\sum^{N}_{i,j=1} a_{i,j}(x) \frac{\partial
f}{\partial x_i}\frac{\partial f}{\partial x_j} +V f^2 \right)dx \
, $$ initially defined on $C^\infty_0(\Omega)$; where
$(a_{i,j}(x))_{n\times n}$ is a measurable symmetric uniformly
elliptic matrix such that
$$\sum^N_{i,j=1} a_{i,j}(x)\xi_{i} \xi_j\ge |\xi|^2$$ and $V$ is a potential on
$\Omega$ such that \be \label{hd1} V(x)=V_1(x)+V_2(x)  ~ , \ee
where \be \label{hd2} |V_1(x)|\le \frac{1}{4d^2(x)} ~ , ~~~
V_2(x)\in L^p(\Omega), ~ p>\frac{N}{2} \ . \ee

We also suppose that \be \label{hd3}\lambda_1:=\inf_{0\neq \varphi
\in C^\infty_0(\Omega)} \frac{Q(\varphi)}{\int_{\Omega} \varphi^2
dx }
>0 , \ee and that to $\lambda_1$ there corresponds a positive
eigenfunction $\varphi_1$ satisfying for all $x\in \Omega$ the
following estimate, \be \label{hd4} c_1\  d^{\frac{\alpha}{2}}(x)
\le \varphi_1(x)\le c_2 \ d^{\frac{\alpha}{2}}(x)\ , ~~ \mbox{ for
some } ~~~\alpha\ge 1\ee and for $c_1, ~ c_2$ two positive
constants.

Thus $\widetilde{E}$ is defined on the closure of
$C^\infty_0(\Omega)$ with respect to the norm defined by the
quadratic form $Q$. Then as before it can be shown that
$\widetilde{E}$ is a well defined nonnegative self-adjoint
operator on $L^2(\Omega)$ such that for every $t>0$,
$e^{-\widetilde{E}t}$ has a integral kernel, that is
$e^{-\widetilde{E}t}u_0(x):=\int_{\Omega}
\widetilde{e}(t,x,y)u_0(y) dy$ and if $N\ge 3$ a Green function
$G_{\widetilde{E}}(x,y)=\int^\infty_0 \widetilde{e}(t,x,y) dt$
denoting the kernel of $\widetilde{E}^{-1}$.

\begin{theorem}
For $N\ge 3$, let $\Omega\subset \R^N$ be a smooth bounded domain.
Suppose that (\ref{hd1}), (\ref{hd2}), (\ref{hd3}) and (\ref{hd4})
are satisfied. Then there exist two positive constants $C_1, C_2$,
with $C_1\le C_2$, such that for any $x,y\in \Omega$
$$C_1 \min\left\{\frac{1}{|x-y|^{N-2}}, \frac{d^{\frac{\alpha}{2}}(x)d^{\frac{\alpha}{2}}(y)}{|x-y|^{N+\alpha-2}}\right\}
\le G_{\widetilde{E}}(x,y) \le C_2
\min\left\{\frac{1}{|x-y|^{N-2}}, \frac{d^{\frac{\alpha}{2}}(x)
d^{\frac{\alpha}{2}}(y)}{|x-y|^{N+\alpha-2}}\right\}.$$
\end{theorem}

Davies conjecture is stated under slightly stronger assumptions on
$V$ than (\ref{hd2}) and on $\varphi_1$ (only when $\alpha=2$).

{\em Proof:} We note that we have $E_1:=-\frac{1}{\varphi_1^2} div
(\varphi_1^2 \nabla)\equiv U^{-1}(\widetilde{E}-\lambda_1)U$,
 $U:L^2(\Omega, \varphi_1^2 dx )\to L^2(\Omega)$ being  the unitary operator
  $Uw:=\varphi_1 w$; hence we have the following relationship between heat
  kernels \be \label{hhh} \widetilde{e}(t,x,y)=\varphi_1(x) \varphi_1(y)e_1(t,x,y)e^{-\lambda_1
  t} \ .\ee Due to (\ref{hd4}) we can apply Theorem
\ref{stab} to the operator $E_1$. Hence due to (\ref{hhh}) for two
positive constants $C_1 \leq C_2$, we have for small time
\be\label{hhh1} C_1 \min\left\{1,
\frac{d^{\frac{\alpha}{2}}(x)d^{\frac{\alpha}{2}}(y)}{t^{\frac{\alpha}{2}}}
\right\} t^{-\frac{N}{2}} e^{-C_2\frac{|x-y|^2}{t}}\le
\widetilde{e}(t,x,y) \le C_2 \min\left\{1,
\frac{d^{\frac{\alpha}{2}}(x)
d^{\frac{\alpha}{2}}(y)}{t^{\frac{\alpha}{2}}}\right\}t^{-\frac{N}{2}}
e^{-C_1\frac{|x-y|^2}{t}}. \ee On the other hand for large time
\be\label{hhh2} C_1 \ d^{\frac{\alpha}{2}}(x)\
d^{\frac{\alpha}{2}}(y) \ e^{-\lambda_1 t} \le
\widetilde{e}(t,x,y)\le
 C_2\  d^{\frac{\alpha}{2}}(x) \ d^{\frac{\alpha}{2}}(y) \ e^{-\lambda_1 t},
\ee for all $x,y\in \Omega$. To obtain this estimate we need to
prove a global Sobolev inequality on $\Omega$, which can be easily
deduced from its local version (\ref{ls}) as well as (\ref{ls2})
with $\lambda=0$ there, by means of a partition of unity as in
\cite{K}. Then the result follows integrating
$\widetilde{e}(t,x,y)$ in the time variable.

\bigskip


\end{document}